\documentclass[11pt,reqno]{amsart}

\usepackage{verbatim}
\usepackage[T1]{fontenc}
\usepackage[utf8]{inputenc}
\usepackage[english]{babel}
\usepackage{csquotes} 
\usepackage{color} 
\usepackage[a4paper,margin=3cm]{geometry}
\usepackage{tikz}  
\usepackage{float}
\allowdisplaybreaks
\usepackage{graphicx}
\usepackage[backend=biber,style=numeric,bibencoding=utf8,language=auto,autolang=other,giveninits=true,doi=false,isbn=false,url=false,maxnames=10,sorting=nty]{biblatex}
\usepackage{amsmath}
\usepackage{amsthm}
\usepackage{amssymb}
\usepackage{mathtools}
\usepackage{esint} 
\usepackage{mathrsfs} 
\usepackage{faktor} 
\usepackage{dsfont} 
\usepackage[dvipsnames]{xcolor}
\usepackage{array}
\usepackage{hhline}
\usepackage[normalem]{ulem}
\usepackage{enumitem} 
\usepackage{caption} 
\usepackage[hyperfootnotes=false]{hyperref} 
\hypersetup{
	colorlinks = true,
	linkcolor = {blue},
	urlcolor = {red},
	citecolor = {blue}
}
\usepackage[nameinlink,capitalise,sort]{cleveref}
\crefname{equation}{}{} 
\crefname{enumi}{}{} 

\crefname{figure}{Figure}{Figures}

\newtheorem{theorem}{Theorem}[section]

\newtheorem{definition}{Definition}[section]
\newtheorem{remark}{Remark}[section]
\newtheorem{proposition}{Proposition}[section]
\newtheorem{corollary}{Corollary}[section]
\newtheorem{example}{Example}[section]

\numberwithin{equation}{section}

\newcommand{\R}{\mathbb R}

\newcommand{\be}{\begin{equation}}
\newcommand{\ee}{\end{equation}}
\newcommand{\ba}{\left(\begin{array}{c}}
\newcommand{\ea}{\end{array}\right)}
\newcommand{\bea}{\begin{eqnarray}}
\newcommand{\eea}{\end{eqnarray}}
\newcommand{\bee}{\begin{eqnarray*}}
\newcommand{\eee}{\end{eqnarray*}}
\newcommand{\ben}{\begin{enumerate}}
\newcommand{\een}{\end{enumerate}}

\begin{document}

\title[KdV on a $\mathcal Y$ junction]{Dynamics of solutions to the Korteweg--de Vries equation on a $\mathcal{Y}$-junction}

\author[M.~Cavalcante]{M\'arcio Cavalcante}
\address[M.~Cavalcante]{Instituto de Matem\'atica, Universidade Federal de Alagoas, 
Av.~Lourival Melo Mota, S/N -- Tabuleiro do Martins, 57072--900 Macei\'o -- AL, Brazil.}
\email{marcio.melo@im.ufal.br} 

\author[N.~De Nitti]{Nicola De Nitti} 
\address[N.~De Nitti]{Politecnico di Bari, Dipartimento di Meccanica, Matematica e Management, Via E.~Orabona 4, 70125 Bari, Italy.}
\email{nicola.denitti@poliba.it}



\begin{abstract}

We study the long-time behavior of solutions to the Korteweg--de Vries equation on a $\mathcal{Y}$-junction, a metric graph consisting of one incoming and two outgoing half-lines joined at a vertex. We first identify vertex conditions of separating-derivative type, together with a balance between the nonlinear coefficients and the coupling parameters, under which the problem is globally well-posed in $L^2(\mathcal{Y})$ for arbitrarily large data. This regime coincides with the one in which the explicit solitary waves of the junction, recently constructed by Mugnolo, Noja and Seifert (2025), exist. Using virial inequalities in the spirit of Mu\~noz and Ponce (2019), we then establish decay, in the sense of $\liminf$, of the $L^2$ norm of the solution in time-dependent regions far from the solitons. In the strongly dominated linear region, which moves away from the vertex, the three edges decouple and the decay holds for every admissible choice of parameters. In the weakly dominated linear region, which stays close to the vertex, the topology of the graph produces boundary contributions at the vertex, and the decay is obtained under an algebraic condition on the coupling parameters that makes these contributions cancel. Finally, we discuss adaptations to general star graphs.

\end{abstract}

\maketitle
\newlist{thmenum}{enumerate}{1}
\setlist[thmenum]{label=(\roman*),ref=\roman*}

\section{Introduction}

The well-known \textit{Korteweg--de Vries (KdV) equation}
\begin{align*}
\partial_t u + \partial_{xxx} u + 2u\,\partial_x u = 0,
\qquad (t,x)\in \mathbb{R}\times\mathbb{R},
\end{align*}
was studied by Boussinesq \cite{Boussinesq1877} in 1877 as a model for long waves propagating on a shallow water surface. Boussinesq's primary aim was to derive an equation with a smooth traveling wave solution, inspired by the wave observed in 1834 by Scott Russell \cite{Russell1845}. 
 Later, it was rediscovered by Korteweg and de Vries \cite{KortewegdeVries1895}.  The KdV equation is now commonly accepted as a mathematical model
for the unidirectional propagation of small-amplitude long waves in nonlinear dispersive
systems. In particular, the KdV equation serves not only as a model for surface water waves but also as a model to study blood pressure waves in large arteries: Shi et al.~\cite[Conclusions and Recommendations]{ShiLawfordHose2011} highlight that while most 1D cardiovascular models neglect vessel wall inertia, some researchers (see \cite{CrepeauSorine2007,Demiray1998,LalegKiratiCrepeauSorine2007,Yomosa1987}) have included this factor, leading to KdV-type dynamics and solitary waves. 

The Cauchy problem for KdV on the real line, the torus, half-lines, and bounded intervals has been extensively studied in recent years (see \cite{Cavalcante2018} and the references therein). 

Progress on the KdV equation on star graphs has been slower than on other domains. The linearized version was studied in the works \cite{SobirovUeckerAkhmedov2015,AkhmedovSobirovEshimbetov2017,SobirovAkhmedovKarpovaJabbarova2015}. Furthermore, Mugnolo, Noja, and Seifert \cite[Theorems~3.9 and~3.16]{MugnoloNojaSeifert2018} obtained a characterization of all boundary conditions under which the Airy-type evolution equation generates a contractive or unitary semigroup of operators. In \cite{Cavalcante2018}, Cavalcante obtained local well-posedness for the Cauchy problem associated with the KdV equation on a metric star graph with three semi-infinite edges—consisting of one negative half-line and two positive half-lines attached to a common vertex—for two classes of boundary conditions.  Parada, Crépeau, and Prieur \cite{ParadaCrepeauPrieur2024} studied the stabilization of the KdV equation in a network with bounded and unbounded edges. In \cite{MugnoloNojaSeifert2025}, the authors studied the existence of solitary waves for the KdV equation on a metric star graph in terms of the coefficients of the equation on each edge, the coupling condition at the central vertex, and the speeds of the traveling waves. Results for the KdV equation on star graphs with bounded edges can be found in the works \cite{CapistranoFilhoParadaSantosdaSilva2025,ParadaCrepeauPrieur2022b,ParadaCrepeauPrieur2022a,Parada2024}. Smith \cite{Smith2026} studied the small-amplitude linearization of the KdV equation on the line with a local defect, where scattering waves are represented by a metric graph domain adjoined at one point.

Finally, linear and nonlinear instability results for stationary solutions on balanced graphs for the KdV equation were obtained in \cite{AnguloPavaCavalcante2021,AnguloPavaCavalcante2025}.

\subsection{Dynamics of solutions to the KdV equation}
To motivate our work, we recall decay results for solutions to the KdV equation on $\mathbb{R}$ and $\mathbb{R}^+$. 

In \cite[Theorem~1.5]{MunozPonce2019}, the authors studied the dynamics of solutions to KdV-type equations on the real line. In particular, they considered the Cauchy problem
\begin{align*}
\begin{cases}
\partial_t u + \partial_{xxx} u + 2u\,\partial_x u = 0,
& (t,x)\in \mathbb R \times \mathbb R,\\
u(0,x)=u_0(x),& x\in \mathbb R.
\end{cases}
\end{align*}
They proved that solutions decay in a sublinear region. More precisely, if
\begin{align*}
u \in C(\mathbb R;H^1(\mathbb R)) \cap L^\infty(\mathbb R;L^1(\mathbb R)),\qquad \sup_{t\in\mathbb R}\|u(t,\cdot)\|_{H^1(\mathbb R)}<\varepsilon,
\end{align*}
where $\varepsilon>0$ is sufficiently small, then
\begin{align*}
\lim_{t\to \pm\infty} \|u(t,\cdot)\|_{H^1(I(t))}=0,
\end{align*}
where the time-dependent interval \(I(t)\) is given by
\begin{align*}
I(t) \coloneqq
\left(
-\frac{C|t|^{1/2}}{\log |t|},
\frac{C|t|^{1/2}}{\log |t|}
\right),
\qquad |t|\geq 2,
\end{align*}
and \(C>0\) is arbitrary.

The ideas developed by Muñoz and Ponce in \cite{MunozPonce2019} and later extended to the Zakharov--Kuznetsov equation in \cite{MendezMunozPobletePozo2021} have been applied to several other models on Euclidean domains. Further results and references are given in \cite[Section~1]{FreireMunoz2025} and \cite{ZhouLinares2026}.

Few results concern the decay of solutions to nonlinear dispersive equations on the half-line. In \cite[Theorem~1.2(2)--(3) and Remark~1.4]{CavalcanteMunoz2023}, decay properties were established for solutions of the initial-boundary value problem on the half-line with homogeneous boundary conditions,
\begin{align*}
\begin{cases}
\partial_t u+\partial_x(\partial_{xx} u+u^2)=0  ,& (t,x)\in (0,\infty)\times\mathbb R_+,\\
u(0,x)=u_0(x),                                   & x\in \mathbb R_+,\\
u(t,0)=0,                                     & t\in(0,\infty).
\end{cases}
\end{align*}
For the solutions in the energy space considered in \cite[Theorem~1.2(2)--(3) and Remark~1.4]{CavalcanteMunoz2023}, the decay results take the form
\begin{itemize}
\item \emph{Weakly dominated linear region:} For any $\delta_0>0$ sufficiently small, 
	\begin{align*}
	\liminf_{t\rightarrow +\infty} \|u(t,\cdot)\|_{H^1(t^{-\delta_0}, t^{4/13-2\delta_0})} =0, \qquad \lim_{t\rightarrow +\infty} \|u(t,\cdot)\|_{H^1(0, t^{-\delta_0} )} =0.
	\end{align*}
	\item \emph{Strongly dominated linear region:} For any $C_1,C_2>0$ and $\frac4{13}\leq n<1$, 
	\begin{align*}
	\liminf_{t\rightarrow +\infty} \|u(t,\cdot)\|_{L^\infty(C_1t^n-C_2,C_1t^n+C_2)}=0.
\end{align*}
\end{itemize}

Recently, in \cite{CamposCardosoCavalcante2026}, the orbital and asymptotic stability of solitons were established for
a class of solutions with non-homogeneous boundary conditions, generalizing some results
of  \cite{CavalcanteMunoz2019} and \cite{CavalcanteMunoz2023}.

Motivated by these results, we study the long-time behavior of solutions to the KdV equation on a $\mathcal{Y}$-junction.

\subsection*{Outline} The paper is organized as follows. \cref{sec:formulation} is devoted to the mathematical formulation of the problem and the statements regarding the existence of solutions for the considered setting. In \cref{sec:results}, we present the main decay results for the KdV equation on the $\mathcal Y$-junction. \crefrange{sec:sublinear}{sec:weak} are devoted to the 
proofs of decay in the strongly and weakly dominated linear regions, respectively. Finally, in \cref{sec:general}, we discuss how the approach developed here can be generalized to other star graphs.

\section{Mathematical Formulation and Well-posedness} 
\label{sec:formulation}

We first specify the boundary conditions at the vertex. Unlike for Schrödinger equations, no canonical choice of boundary conditions for KdV equations on star graphs is generally accepted. Mugnolo, Noja, and Seifert \cite[Theorems~3.9 and~3.16]{MugnoloNojaSeifert2018} characterize the conditions under which the Airy operator generates a contractive or unitary semigroup. A more specific subclass of interest consists of separating derivative boundary conditions. In this direction, existence results for $\mathcal{Y}$-junctions and general star graphs have been established in \cite{Cavalcante2018} and \cite{CavalcanteMarques2025}, respectively.

In this work,   we consider real-valued solutions of the KdV equation on a junction $\mathcal{Y}$ consisting of three semi-infinite edges: one negative half-line and two distinct positive half-lines whose endpoints are identified at a common vertex. More precisely,
\begin{align}\label{KDVY}
\begin{cases}
\partial_t u + \partial_{xxx} u + \gamma_1\, u\,\partial_x u = 0,
& (t,x)\in(0,T)\times(-\infty,0),\\[4pt]
\partial_t v + \partial_{xxx} v + \gamma_2v\,\partial_x v = 0,
& (t,x)\in(0,T)\times(0,+\infty),\\[4pt]
\partial_t w + \partial_{xxx} w + \gamma_3 w\,\partial_x w = 0,
& (t,x)\in(0,T)\times(0,+\infty).
\end{cases}
\end{align}
We impose the initial conditions
\begin{align}\label{initial}
u(0,x)=u_0(x),\qquad v(0,x)=v_0(x),\qquad w(0,x)=w_0(x),
\end{align}
where
\begin{align*}
(u_0,v_0,w_0)\in H^{s}(\R^-)\times H^{s}(\R^+)\times H^{s}(\R^+)\eqqcolon H^{s}(\mathcal{Y}).
\end{align*}
Following \cite[Corollary~1.1]{Cavalcante2018}, we impose the boundary conditions
\begin{align}\label{cond}
\begin{cases}
u(t,0)=\alpha_2v(t,0)=\alpha_3w(t,0),& t\in (0,T),\\
\partial_x u(t,0)=\beta_2\,\partial_x v(t,0)+\beta_3\,\partial_x w(t,0),& t\in(0,T),\\
\partial_{xx} u(t,0)=\frac{1}{\alpha_2}\,\partial_{xx} v(t,0)+\frac{1}{\alpha_3}\,\partial_{xx} w(t,0), & {t\in(0,T),}
\end{cases}
\end{align}
where $\alpha_2,\alpha_3\in\mathbb R\setminus\{0\}$, $\beta_2,\beta_3\in\mathbb R$, and $\gamma_1,\gamma_2,\gamma_3\in\mathbb R$ are fixed constants.
\begin{figure}[h]
    \centering
    \begin{tikzpicture}[line cap=round]
        \draw[dashed, black!85, very thick] (-6,0) -- (-4,0);
        \draw[thick, black!85, very thick] (-4,0) -- (0,0);
        \draw[black!85, ->, thick] (-4.6,0) -- (-4.5,0);
        \draw[black!85, ->, thick] (-3,0) -- (-2.5,0);
        \draw[thick, black!55, very thick] (0,0) -- (30:3);
        \draw[dashed, black!55, very thick] (30:3) -- (30:4);
        \draw[black!55, ->, thick] (30:1) -- (30:1.5);
        \draw[black!55, ->, thick] (30:3.5) -- (30:3.75);
        \draw[thick, black!55, very thick] (0,0) -- (-30:3);
        \draw[dashed, black!55, very thick] (-30:3) -- (-30:4);
        \draw[black!55, ->, thick] (-30:1) -- (-30:1.5);
        \draw[black!55, ->, thick] (-30:3.5) -- (-30:3.75);
        \node at (-4.5,0.5) {$( -\infty, 0)$};
        \node at (4.5,1.7) {$(0, +\infty)$};
        \node at (4.5,-1.7) {$(0, +\infty)$};
        \filldraw (0,0) circle (2pt);

    \end{tikzpicture}
    \caption{$\mathcal{Y}$-junction}
    \label{fig-y-junction}
\end{figure}
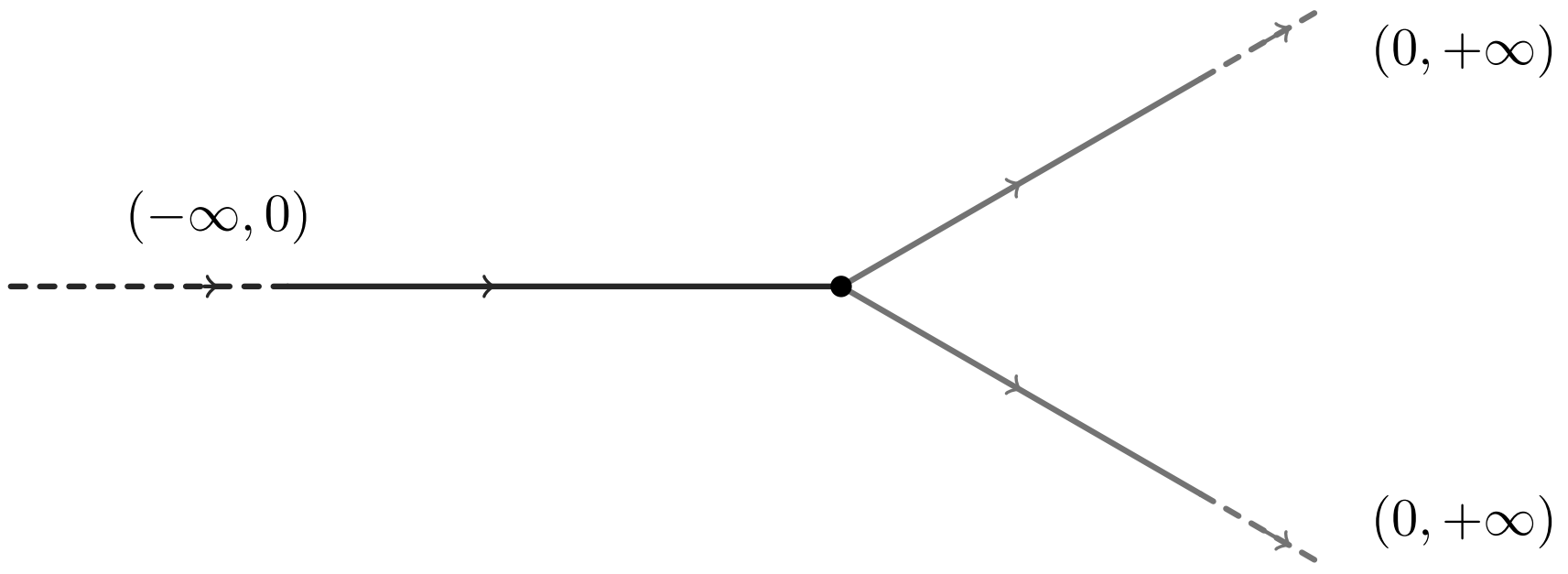

We write the system \cref{KDVY} in edge notation as
\begin{align*}
 \begin{cases}
 \partial_t\mathcal U_e+\partial_{xxx}\mathcal U_e+\gamma_e\mathcal U_e\partial_x\mathcal U_e=0,& (t,x)\in(0,T)\times e,\\
 \mathcal U_e(0,x)=\mathcal U_{e,0}(x),&x\in e,
 \end{cases}
\end{align*}
where $e$ ranges over the three edges, labelled $1,2,3$, and $\mathcal U_1=u$, $\mathcal U_2=v$, and $\mathcal U_3=w$. We write $\mathcal U=(u,v,w)$ for the full solution. An integral over a region of $\mathcal Y$ denotes the sum of the integrals over its three edge components.

The local existence result in \cite[Corollary~1.1]{Cavalcante2018} follows from the boundary forcing method for the KdV equation on half-lines developed in \cite{CollianderKenig2002,Holmer2006}. Although \cite{Cavalcante2018} treats unit nonlinear coefficients, its contraction argument in \cite[Section~5]{Cavalcante2018} also applies to fixed real coefficients $\gamma_1,\gamma_2,\gamma_3$: these coefficients multiply the bilinear terms and do not change the boundary matrix. We use this extension in \cref{local}.

\begin{proposition}\label{local} Let $-\frac{1}{2}<s<\frac{3}{2}$ with $s\neq\frac{1}{2}$, let $\gamma_1,\gamma_2,\gamma_3\in\mathbb R$, and let $\alpha_2,\alpha_3\in\mathbb R\setminus\{0\}$ and $\beta_2,\beta_3\in\mathbb R$ satisfy
\begin{align*}
\frac{1}{\alpha_2^2}+\frac{1}{\alpha_3^2}+\frac{\beta_3}{\alpha_3}+\frac{\beta_2}{\alpha_2}\neq-1.
\end{align*}
Assume that $(u_0,v_0,w_0)\in H^s(\mathcal Y)$. If $\frac12<s<\frac32$, assume also that $u_0(0)=\alpha_2v_0(0)=\alpha_3w_0(0)$. Then there exists a positive time $T$ and a distributional solution $(u,v,w)\in C([0,T];H^s(\mathcal Y))$ of \cref{KDVY,initial,cond}. The data-to-solution map is locally Lipschitz continuous on the set of compatible data, with values in $C([0,T];H^s(\mathcal Y))$. The solution is unique in the auxiliary solution class used in \cite[Section~5]{Cavalcante2018}.
\end{proposition}
In \cite[Theorem~2.1]{CavalcanteMarques2025}, the existence of solutions for the KdV equation was studied on more general star graphs with an arbitrary number of edges. In the recent paper \cite{MugnoloNojaSeifert2025}, the KdV equation on a metric star graph was studied, and the existence of solitary waves on the metric graph was investigated in terms of the coefficients of the equation on each edge, the coupling condition at the central vertex of the star, and the speeds of the traveling wave. In particular, discontinuous solitons were obtained in the case of a $\mathcal Y$ junction.

The existence of solitons raises the question of which regions of the $\mathcal Y$-junction allow decay of solutions to \cref{KDVY} as $t\to\infty$.  To our knowledge, the asymptotic behavior of solutions to the KdV equation on a $\mathcal Y$-junction has not been studied before.

We first identify conditions on the parameters under which the local solutions can be extended globally in time.

\begin{proposition}\label{zamp1}
Let $\alpha_2,\alpha_3,\beta_2,\beta_3\in\mathbb{R}\setminus\{0\}$ satisfy
the condition of \cref{local}, and let
\begin{align*}
\gamma_1=\frac{1}{\alpha_2^{3}}+\frac{1}{\alpha_3^{3}},\qquad
\gamma_2=\gamma_3=1,\qquad
\beta_2^2+\beta_3^2\le 1 .
\end{align*}
Then the local solution of the Cauchy problem
\cref{KDVY,initial,cond} in the space $C([0,T];L^2(\mathcal{Y}))$,
obtained in \cref{local}, can be extended to any time $T>0$.
\end{proposition}

\begin{proof}
    We obtain the mass identity \cref{tec_cubic} formally by multiplying each component equation by the corresponding component of the solution and integrating by parts on each edge:
 \begin{align}\label{tec_cubic}
\begin{aligned}
&\frac12 \frac{\mathrm{d}}{\mathrm{d}t}\left(
\int_{-\infty}^{0} u^2(t,x)\,\mathrm{d}x
+ \int_0^{\infty} v^2(t,x)\,\mathrm{d}x
+ \int_0^{\infty} w^2(t,x)\,\mathrm{d}x
\right)\\
&= -u(t,0)\,\partial_{xx} u(t,0)
+ v(t,0)\,\partial_{xx} v(t,0)
+ w(t,0)\,\partial_{xx} w(t,0)\\
&\quad + \frac12\left(
(\partial_x u(t,0))^2
- (\partial_x v(t,0))^2
- (\partial_x w(t,0))^2
\right)\\
&\quad - \frac13\left(
\gamma_1 u^3(t,0) - v^3(t,0) - w^3(t,0)
\right)\\
&\eqqcolon A + B + C.
\end{aligned}
\end{align}
 Using \cref{cond}, the coefficient relation, and the Cauchy--Schwarz inequality, we obtain $A=C=0$ and $B\leq0$.

We justify the identity at $L^2$ regularity using \cite[Section~5, Eq.~(5.23)]{Cavalcante2018}. The fixed-point construction admits approximations by finite Picard iterates with compatible $H^1$ data, converging in its $s=0$ solution norm. The trace estimates then give spatially continuous traces of each component and its first two derivatives in $H^{1/3}(0,T)$, $L^2(0,T)$, and $H^{-1/3}(0,T)$, respectively. The first two trace norms vanish at spatial infinity, by approximation from the $H^1$ iterates and Sobolev interpolation. On bounded edge segments, local smoothing and the one-dimensional Sobolev inequality give $L^4$ integrability of each component and $L^2$ integrability of its derivative. Spatial mollification therefore justifies \cref{tec_cubic}: the nonlinear commutator tends to zero by the Cauchy--Schwarz inequality, and the mixed boundary products converge by $H^{1/3}$--$H^{-1/3}$ duality. The fluxes at infinity vanish, and \cref{cond} cancels the mixed and cubic vertex terms. Thus
\begin{align*}
\|\mathcal U(t,\cdot)\|_{L^2(\mathcal Y)}^2
&+\int_0^t\Big[(\partial_x v(s,0))^2+(\partial_x w(s,0))^2
-(\beta_2\partial_x v(s,0)+\beta_3\partial_x w(s,0))^2\Big]\,\mathrm{d}s\\
&=\|(u_0,v_0,w_0)\|_{L^2(\mathcal Y)}^2.
\end{align*}
The integrand is non-negative by the Cauchy--Schwarz inequality. Since the local existence time in \cite[Section~5]{Cavalcante2018} is controlled by the $L^2$ norm, iteration extends the solution to every finite time.
\end{proof}

\begin{remark}\label{rk:1}
The condition $\gamma_1 = \alpha_2^{-3} + \alpha_3^{-3}$ cancels the cubic boundary flux in \cref{tec_cubic}, allowing the $L^2$-mass bound in \cref{zamp1}.
\end{remark}

\begin{remark}[$L^2(\mathcal Y)$ and $H^1(\mathcal Y)$ bounds]
The cancellation in \cref{rk:1} controls the $L^2$-norm. It does not provide a conserved quantity controlling the $H^1$-norm under \cref{cond}. The decay results in \cref{theo_1,teo_2} use only the $L^2$-bound.
\end{remark}

The proof of \cref{zamp1} also gives the mass dissipation bound in \cref{cor_1}.
\begin{corollary}\label{cor_1} Let  $\gamma_1=\frac{1}{\alpha_2^3}+\frac{1}{\alpha_3^3}$, $\gamma_2=\gamma_3=1$, and $\beta_2^2 + \beta_3^2 \le 1$. 
 Then solutions of the Cauchy problem \crefrange{KDVY}{cond} satisfy
\begin{align*}
 &\int_{-\infty}^{0} u^2(t,x)\,\mathrm{d}x
+ \int_0^{\infty} v^2(t,x)\,\mathrm{d}x
+ \int_0^{\infty} w^2(t,x)\,\mathrm{d}x\\ 
&\leq \int_{-\infty}^{0} u_0^2\,\mathrm{d}x
+ \int_0^{\infty} v_0^2\,\mathrm{d}x
+ \int_0^{\infty} w_0^2\,\mathrm{d}x.
\end{align*}
\end{corollary}

We review the solitary waves constructed in \cite{MugnoloNojaSeifert2025} to identify regions in which decay cannot hold for all solutions to \cref{KDVY}. We first define travelling and solitary waves on the $\mathcal{Y}$-junction.

\begin{definition}\label{def:travelling-wave}
We call $\mathcal Q=(Q_1,Q_2,Q_3)$ a \emph{travelling wave} for \cref{KDVY,cond} if there exist profiles $\varphi_i\in C^3(\mathbb R)$ and speeds $c_i\in\mathbb R\setminus\{0\}$ such that $Q_i(t,x)=\varphi_i(x-c_it)$ on the $i$th edge and $\mathcal Q$ satisfies the equations and coupling conditions. We call it a \emph{solitary wave} if each profile is non-constant, has at most one local extremum, and satisfies $\lim_{y\to\pm\infty}\varphi_i(y)=0$.
\end{definition}

The profiles in \cref{def:travelling-wave} may have different traces at the vertex.

The decaying profiles of solitary waves considered in \cite[Theorem~2.7]{MugnoloNojaSeifert2025} have the form
\begin{align}\label{eq:solitary-profile}
\varphi_i(y)\coloneqq\frac{3c_i}{\gamma_i}\operatorname{sech}^2\left(\frac{\sqrt{c_i}}2(y-y_{0,i})\right),
\end{align}
where $c_i>0$, $\gamma_i\neq0$, and $y_{0,i}\in\mathbb R$ for $i\in\{1,2,3\}$. For the coupling conditions \cref{cond}, we obtain the criterion in \cref{Teo_solitons} by comparing the vertex traces of these profiles. We refer to \cite[Proposition~4.5]{MugnoloNojaSeifert2025} for the corresponding construction on star graphs.

\begin{proposition}[Existence of solitons on a $\mathcal Y$-junction]\label{Teo_solitons}
Let $\alpha_2=\alpha_3=1/a$, where $a\in\mathbb R\setminus\{0\}$, let $\gamma_1,\gamma_2,\gamma_3\in\mathbb R\setminus\{0\}$, and assume $\beta_2^2+\beta_3^2\leq1$. The profiles in \cref{eq:solitary-profile} define a solitary wave for \cref{KDVY,cond} if and only if
\begin{align*}
\frac{\beta_2}{c_2}+\frac{\beta_3}{c_3}&=\frac1{ac_1},
& c_1,c_j&>0,\\
\frac{y_{0,1}}{c_1}&=\frac{y_{0,j}}{c_j},
& \sqrt{c_1}\,c_1&=\sqrt{c_j}\,c_j,\\
a\frac{c_1}{\gamma_1}&=\frac{c_j}{\gamma_j},
& \frac1{c_1^2}&=a^2\left(\frac1{c_2^2}+\frac1{c_3^2}\right),
\end{align*}
where the relations involving $j$ hold for $j\in\{2,3\}$.
\end{proposition}
The factor $1/a$ in the derivative condition follows by differentiating $u(t,0)=v(t,0)/a=w(t,0)/a$ in time and using $\partial_t Q_i=-c_i\partial_x Q_i$.
\begin{remark}[Rigidity of solitary waves on the $\mathcal{Y}$-junction]\label{rem:rigidity}
The conditions in \cref{Teo_solitons}, together with $\gamma_2=\gamma_3=1$, imply
\begin{align*}
c_1=c_2=c_3\eqqcolon c>0,\qquad a^2=\frac12,\qquad\gamma_1=a,
\end{align*}
by the injectivity of $c\mapsto c^{3/2}$ on $(0,\infty)$ and the last two relations in that proposition. Consequently,
\begin{align*}
\alpha_2=\alpha_3=\pm\sqrt2,\qquad\gamma_1=\pm\frac1{\sqrt2},
\end{align*}
with matching signs, and in both cases $\gamma_{1}=\alpha_{2}^{-3}+\alpha_{3}^{-3}$, so that the
coefficient condition in \cref{zamp1} is satisfied. Thus \cref{ex_1}, up to the sign of $a$, the speed $c$, and a common translation, describes all solitary waves in this class.
\end{remark}

\begin{example}\label{ex_1} Let $\alpha_2 = \alpha_3 = \sqrt{2}$, $\gamma_1 = \frac{1}{\sqrt{2}}$, $\gamma_2 = \gamma_3 = 1$, so that, by \cref{rem:rigidity}, the wave speeds are necessarily
identical across all edges, $c_{1}=c_{2}=c_{3}\eqqcolon c>0$.

The derivative condition at the vertex requires $\frac{\beta_2}{c_2} + \frac{\beta_3}{c_3} = \frac{\alpha_2}{c_1}$, which reduces to:
\begin{align*}
    \beta_2 + \beta_3 = \sqrt{2}.
\end{align*}
The mass inequality in \cref{zamp1} also requires
\begin{align*}
    \beta_2^2 + \beta_3^2 \le 1.
\end{align*}

By the Cauchy--Schwarz inequality, the minimum of $\beta_2^2 + \beta_3^2$ subject to $\beta_2 + \beta_3 = \sqrt{2}$ is $1$, attained only when $\beta_2 = \beta_3$. Therefore, the only derivative coupling compatible with these conditions is
\begin{align*}
    \beta_2 = \beta_3 = \frac{\sqrt{2}}{2}.
\end{align*}

For a common translation parameter $y_0$, the solitary wave is
\begin{align*}
    \begin{cases}
        u(t,x) = 3\sqrt{2}c \operatorname{sech}^2\left(\frac{\sqrt{c}}{2}(x - ct - y_0)\right), & x \in (-\infty, 0), \\
        v(t,x) = 3c \operatorname{sech}^2\left(\frac{\sqrt{c}}{2}(x - ct - y_0)\right), & x \in (0, \infty), \\
        w(t,x) = 3c \operatorname{sech}^2\left(\frac{\sqrt{c}}{2}(x - ct - y_0)\right), & x \in (0, \infty).
    \end{cases}
\end{align*}
This solution satisfies the boundary conditions \cref{cond}.

\begin{figure}[ht] 
        \centering      \includegraphics[width=0.5\linewidth]{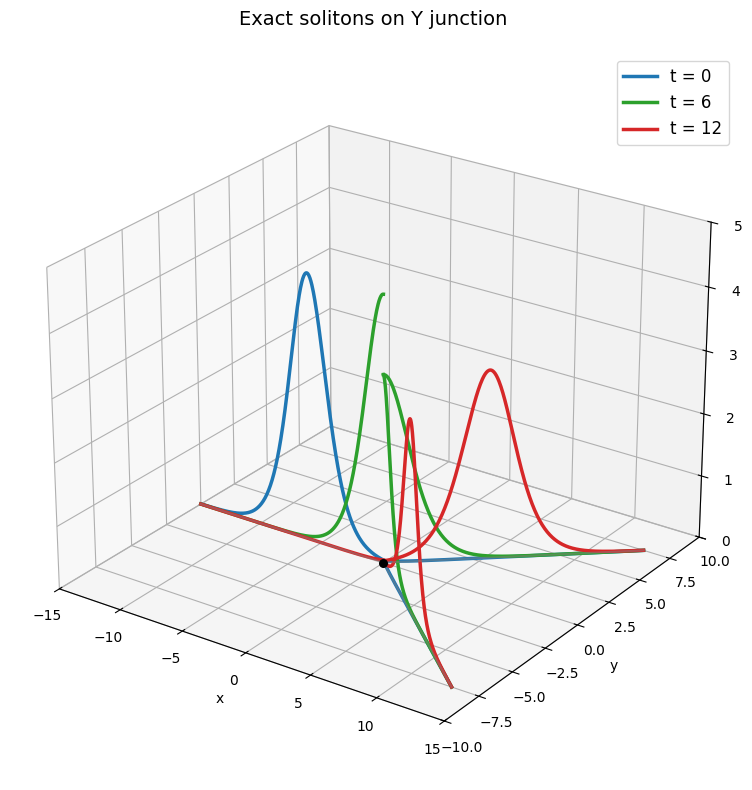}
        \caption{3D visualization of the symmetric soliton propagation at three different time steps ($t=0, 6, 12$), showing the wave traversing the vertex and splitting into the outgoing edges.}
    \label{fig2:soliton_3d}
\end{figure}

\end{example}

We study decay in time-dependent regions of the $\mathcal Y$-junction that do not follow a solitary wave. We establish decay in the sense of a $\liminf$ as $t\to\infty$ in the weakly and strongly dominated linear regions.

\section{Main results}
\label{sec:results}

In this section, we state the main results, which are divided into two parts.

Bona and Cascaval \cite[Section~4 and Figure~4]{BonaCascaval2008} numerically studied nonlinear dispersive waves on networks for the Benjamin--Bona--Mahony (BBM) equation. In their numerical simulations of a $\mathcal Y$-junction, they observed that an incoming solitary wave scatters at the bifurcation, generating a transmitted wave alongside a negative reflected wave. Mugnolo et al. \cite{MugnoloNojaSeifert2025} constructed explicit travelling waves for the KdV equation on graphs. We study local $L^2$-decay for solutions satisfying the assumptions in \cref{theo_1,teo_2}.

\begin{figure}[ht]
\centering
\begin{tikzpicture}[x=1cm,y=1cm,font=\small,every path/.style={draw=.},every node/.style={text=.}]
    \node[anchor=west] at (0,0.85) {Strongly dominated region $\Omega(t)$};
    \draw[->] (0,0) -- (10.8,0) node[right] {$r=|x|$};
    \fill[fill=.] (0,0) circle (0.045);
    \node[below] at (0,-0.1) {$0$};
    \draw[line width=2pt] (3.8,0) -- (9.2,0);
    \draw[fill=white] (3.8,0) circle (0.06);
    \draw[fill=white] (9.2,0) circle (0.06);
    \draw (6.5,-0.1) -- (6.5,0.1);
    \node[above] at (6.5,0.12) {$At^n$};
    \node[below] at (3.8,-0.12) {$At^n-\dfrac{Bt^b}{\log t}$};
    \node[below] at (9.2,-0.12) {$At^n+\dfrac{Bt^b}{\log t}$};

    \node[anchor=west] at (0,-1.45) {Weakly dominated region $\Omega_0(t)$};
    \draw[->] (0,-2.3) -- (10.8,-2.3) node[right] {$r=|x|$};
    \fill[fill=.] (0,-2.3) circle (0.045);
    \node[below] at (0,-2.4) {$0$};
    \draw[line width=2pt] (1.5,-2.3) -- (9.2,-2.3);
    \draw[fill=white] (1.5,-2.3) circle (0.06);
    \draw[fill=white] (9.2,-2.3) circle (0.06);
    \node[below] at (1.5,-2.42) {$t^{-\varepsilon_1}$};
    \node[below] at (9.2,-2.42) {$t^{4/13-\varepsilon_2}$};
\end{tikzpicture}
\caption{Decay intervals on each edge at a fixed sufficiently large time $t$, expressed in the coordinate $r=|x|$.}
\label{fig:decay-intervals}
\end{figure}
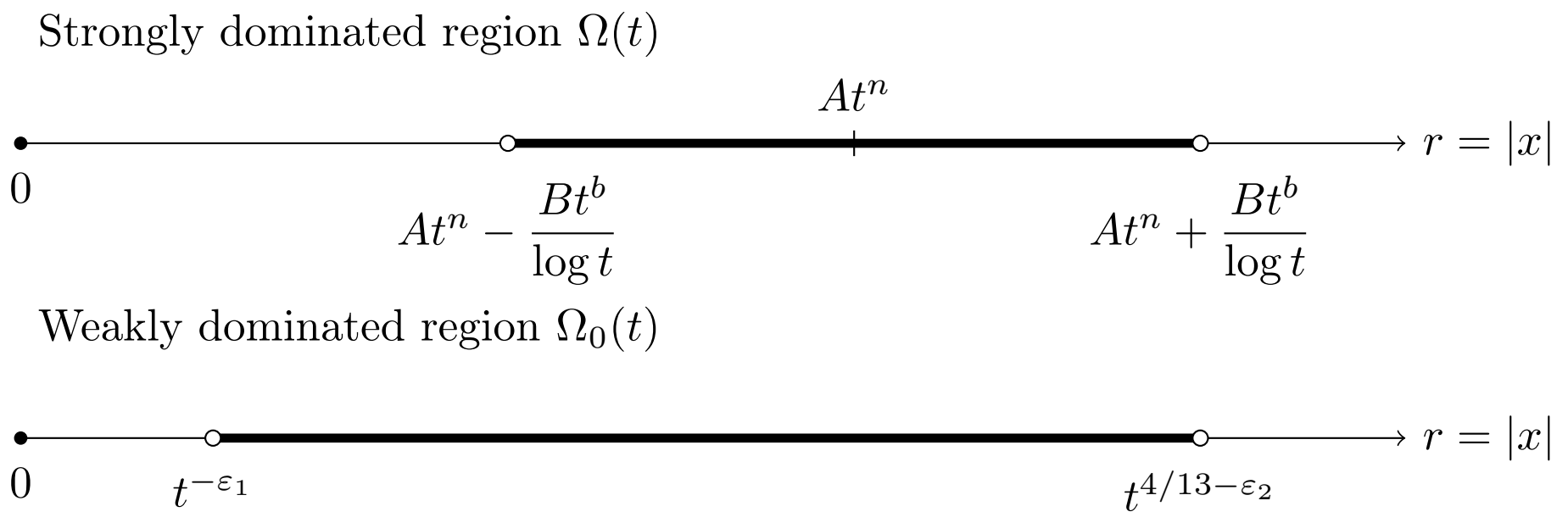

To study the strongly dominated linear region, we define the time-dependent region $\Omega(t)\subset\mathcal{Y}$ illustrated in \cref{fig:decay-intervals}, using the distance $r=|x|$ from the vertex on each edge. More precisely, for sufficiently large $t$, we set
\begin{align*}
\Omega(t)\coloneqq\big(I(t),J(t),J(t)\big),
\end{align*}
where
\begin{align*}
I(t)\coloneqq\Big(-At^{n}-\frac{Bt^{b}}{\log t},\,-At^{n}+\frac{Bt^{b}}{\log t}\Big),
\qquad
J(t)\coloneqq\Big(At^{n}-\frac{Bt^{b}}{\log t},\,At^{n}+\frac{Bt^{b}}{\log t}\Big),
\end{align*}
where $A,B>0$, $0<b<\frac23$, and $b<n\le 1-\frac{b}{2}$ are fixed.

\begin{theorem}[Strongly dominated linear region]\label{theo_1} Let $\alpha_2,\alpha_3,\beta_2,\beta_3\in\mathbb R\setminus\{0\}$ satisfy
\begin{align*}
&\alpha_2^{-2}+\alpha_3^{-2}+\beta_2/\alpha_2+\beta_3/\alpha_3\neq-1,\\
&\gamma_1=\alpha_2^{-3}+\alpha_3^{-3}\neq0,\qquad\gamma_2=\gamma_3=1,\qquad\beta_2^2+\beta_3^2\le1.
\end{align*}
 Let $\mathcal U=(u,v,w)\in C([0,\infty);L^2(\mathcal Y))$ be the global solution of \cref{KDVY,initial,cond} given by \cref{zamp1}.
Assume that $0<b<\frac23$ and $b< n \leq 1-\frac{b}2$. Then, for any fixed $A,B >0$,
\begin{align}\label{Auxiliar_2_bis}
 \liminf_{t \rightarrow+\infty}\int_{\Omega(t)}\mathcal U_{e}^2 (t,x) \,\mathrm{d}x=0.
\end{align}
\end{theorem}
 \medskip
\begin{remark}
The assumption $\gamma_{1}\neq 0$ is only used on the negative edge, where
the coercive term in the virial identity is proportional to $|\gamma_{1}|$
(see \cref{captal-xi_left}). When $\gamma_{1}=0$ the equation on
$(-\infty,0)$ is the linear Airy equation and \cref{Auxiliar_2_bis}
still holds on the two positive edges. The argument does not cover the negative edge in this case.
\end{remark}

\medskip
 We next define the time-dependent weakly dominated linear region $\Omega_0(t)\subset\mathcal{Y}$, also illustrated in \cref{fig:decay-intervals}. Given $0<\varepsilon_1<\varepsilon_2$, we set $\Omega_0(t)  \coloneqq  (I_0(t), J_0(t), J_0(t))$, where $I_0(t) = (-t^{\frac{4}{13}-\varepsilon_2},-t^{-\varepsilon_1})$ and $J_0(t) = (t^{-\varepsilon_1},t^{\frac{4}{13}-\varepsilon_2})$. The distinct endpoint exponents permit the choice of an intermediate exponent in \cref{eq:params}, which is used in \cref{eq:weight-below}.

\begin{theorem}[Weakly dominated linear region]\label{teo_2} Let $\alpha_2,\alpha_3,\beta_2,\beta_3\in\mathbb R\setminus\{0\}$ satisfy
\begin{align*}
&\alpha_2^{-2}+\alpha_3^{-2}+\beta_2/\alpha_2+\beta_3/\alpha_3\neq-1,\\
&\gamma_1=\alpha_2^{-3}+\alpha_3^{-3}\neq0,\qquad\gamma_2=\gamma_3=1,\qquad\beta_2^2+\beta_3^2\le1.
\end{align*}
Assume also that there exist
$A_1,A_2,A_3>0$ such that
\begin{align}\label{eq:positive-multipliers}
    \operatorname{sgn}(\gamma_1)A_1
    +\frac{A_2}{\alpha_2}
    +\frac{A_3}{\alpha_3}=0.
\end{align} Let $\mathcal U=(u,v,w)\in C([0,\infty);L^2(\mathcal Y))$ be the global solution of \cref{KDVY,initial,cond} given by \cref{zamp1}. Then, for every
$0<\varepsilon_1<\varepsilon_2<4/39$,  we have
\begin{align*}
 \liminf_{t \rightarrow+\infty}\int_{\Omega_0(t)}\mathcal U_{e}^2 (t,x) \,\mathrm{d}x=0.
\end{align*}

\end{theorem}
\begin{remark}[On condition \cref{eq:positive-multipliers}]\label{rem:boundary-cancellation}
Unlike the Cauchy problem on the real line, the topology of the
$\mathcal{Y}$-junction produces boundary contributions at the vertex $x=0$ in the
weighted virial functionals of \cref{sec:weak}. On each edge, the sign of the
weight is dictated by the coercivity of the nonlinear term, which in turn fixes the trace
term at the vertex: these are $2\operatorname{sgn}(\gamma_1)\,u(t,0)$, $2v(t,0)$, and $2w(t,0)$
(up to the common factor $\lambda_1^{-2}(t)\eta^{-1}(t)$). The estimates used in the proof do not establish time integrability of these terms separately. We choose the positive multipliers
$A_1,A_2,A_3$ in the combined functional \cref{eq:Xi-total-new}. The algebraic condition
\cref{eq:boundary-cancellation}, together with the transmission condition
$u(t,0)=\alpha_2v(t,0)=\alpha_3w(t,0)$, makes the three traces cancel.
The positivity of the $A_j$ is needed to add the edge inequalities in \cref{eq:total-virial} and to deduce \cref{eq:integrated}.

Since $\gamma_1=\alpha_2^{-3}+\alpha_3^{-3}$ under the assumptions of
  \cref{zamp1}, the three numbers $\operatorname{sgn}(\gamma_1)$,
$\alpha_2^{-1}$, and $\alpha_3^{-1}$ share a common sign whenever $\alpha_2\alpha_3>0$, and
\cref{eq:boundary-cancellation} has no solution with $A_1,A_2,A_3>0$. Conversely, if
$\alpha_2\alpha_3<0$ and $\gamma_1\neq0$, then \cref{eq:boundary-cancellation} always admits positive
solutions, and by homogeneity we may normalize $A_2=1$. Under the other assumptions of \cref{teo_2}, including $\gamma_1\neq0$, condition \cref{eq:positive-multipliers} is therefore equivalent to $\alpha_2\alpha_3<0$. In that case the first condition in \cref{cond} reverses the sign on one outgoing edge. For instance,
$\alpha_2=1$, $\alpha_3=-2$, $\beta_2=\beta_3=1/\sqrt2$ gives $\gamma_1=7/8$ and
$(A_1,A_2,A_3)=(1,1,4)$. In particular, the symmetric configuration of
\cref{ex_1} is not covered by \cref{teo_2}. Whether
$L^2$-decay in $\Omega_0(t)$ holds in that case, where the vertex traces do not cancel,
remains an open question.
\end{remark}

\section{Analysis in the strongly dominated linear region}
\label{sec:sublinear}

We prove \cref{theo_1}, which establishes $L^2$-decay in the sense of $\liminf$ in the region $\Omega(t)$. The centres of its component intervals have distance $At^n$ from the vertex. 

\begin{proof}
\uline{Step 1:} \emph{Virial functionals.} We adapt the proof of \cite[Lemma~5.2]{CavalcanteMunoz2023} to the junction $\mathcal{Y}$, using weights supported away from the vertex to avoid boundary traces.

We choose a non-decreasing function $\psi\in C^\infty(\mathbb R)$ such that $\psi(s)=1$ for $s\ge1$ and $\psi(s)=0$ for $s\le-1$, and set $\varphi\coloneqq\psi'\ge0$.
The function $\varphi$ has compact support in $[-1,1]$. We also choose $\psi$ so that $\psi^{\prime}\geq c_0$ on $[-\tfrac12,\tfrac12]$, where $c_0>0$ is fixed, to obtain a positive lower bound for the weight in the region of interest.

We write $\mathcal U=(u,v,w)$ for the solution of \cref{KDVY,initial,cond}.

On the second edge we consider the virial functional
\begin{align}\label{captal-xi}
\Xi[v](t) \coloneqq \frac{1}{\eta(t)}\int_{0}^\infty v(t,x)\psi\left(\frac{x - \tilde\rho(t)}{\lambda_{1}(t)}\right)\varphi\left(\frac{x- \tilde\rho(t) }{\lambda_{1}^q(t)}\right)\,\mathrm{d}x, 
\end{align}
for $\lambda_1(t)= 4B t^b /\log t $ and $\eta(t)=t^{m}\ln ^{2} t.$ We consider the translation parameter $\tilde\rho(t)$ of the form
\begin{align}\label{rho_new}
\tilde\rho(t) = At^n, \quad \text{with  $A>0$  fixed.}
\end{align}

By the support of $\varphi$, the integral in \cref{captal-xi} is concentrated on the set \begin{align*}\left[At^n-\left(\frac{4Bt^b}{{\log}t}\right)^{q},At^n+\left(\frac{4Bt^b}{{\log}t}\right)^{q}\right].\end{align*}

We choose the exponents so that
\begin{align}\label{cc1_new}
\begin{aligned}
&m,b>0, \qquad m=1-b, \\
&b < n \leq 1-\frac{b}2,
\end{aligned}
\end{align}
and, for some $q_0>0$ small enough,
\begin{align}\label{cc2_new}
b \leq  \frac{2}{2+q}, \quad q=1+q_0, \quad 0<q_0<\frac{n-b}{b}.
\end{align}

We fix $T_0>10$ sufficiently large that $At^n>\lambda_1^q(t)$, $\lambda_1(t)\geq1$, and $\lambda_1^{1-q}(t)\leq\tfrac12$ for $t\geq T_0$. Thus the weights on all three edges are supported away from the vertex. All identities in this proof are understood for $t\geq T_0$. They hold in distributions in time by testing \cref{KDVY} against these smooth weights. Since the solution belongs to $C(\mathbb R^+;L^2(\mathcal Y))$, the resulting right-hand sides are locally integrable, and the virial functionals are locally absolutely continuous.

The Cauchy--Schwarz inequality and \cref{cc2_new} yield
\begin{align*}
\sup_{t\gg 1} |\Xi[v](t)| \lesssim \sup_{t\gg 1} \frac{\left(\lambda_{1}^q(t)\right)^{1/2}}{\eta(t)} \lesssim  \sup_{t\gg 1} \frac1{t^{1-b-\frac12 bq} \ln^{2+\frac{q}2} t} <\infty.
\end{align*}
We set $\tilde x\coloneqq x-\tilde\rho(t)$. We have
\begin{align*}
\begin{aligned}
 \frac{\mathrm{d}}{\mathrm{d}t}\Xi[v](t)&= \frac{1}{\eta(t)}\int_{0}^\infty\partial_t\left(v\psi\left(\frac{\tilde x }{\lambda_{1}(t)}\right)\varphi\left(\frac{\tilde x }{\lambda_{1}^q(t)}\right)\right)\,\mathrm{d}x \\ & \quad -\frac{\eta'(t)}{\eta^{2}(t)}\int_{0}^\infty v\psi\left(\frac{\tilde{x}}{\lambda_{1}(t)}\right)\varphi\left(\frac{\tilde x }{\lambda_{1}^q(t)}\right)\,\mathrm{d}x\\
 &\eqqcolon \Xi_{1}[v](t)+ \Xi_{2}[v](t).
\end{aligned}
\end{align*}
Differentiating the weights gives
\begin{align*}
\begin{aligned}
\Xi_{1}[v](t)&=\frac{1}{\eta(t)}\int_{0}^\infty\partial_t v\psi\left(\frac{\tilde x }{\lambda_{1}(t)}\right)\varphi\left(\frac{\tilde x }{\lambda_{1}^q(t)}\right)\,\mathrm{d}x  \\ & \quad -\frac{\lambda_{1}'(t)}{\lambda_{1}(t)\eta(t)}\int_{0}^\infty v \psi'\left(\frac{\tilde x }{\lambda_{1}(t)}\right)\left(\frac{\tilde x }{\lambda_{1}(t)}\right)\varphi\left(\frac{\tilde x }{\lambda_{1}^q(t)}\right)\,\mathrm{d}x\\
&\quad -q\frac{\lambda_{1}'(t)}{\lambda_{1}(t)\eta(t)}\int_{0}^\infty v \psi\left(\frac{\tilde x }{\lambda_{1}(t)}\right)\left(\frac{\tilde x }{\lambda_{1}^q(t)}\right)\varphi'\left(\frac{\tilde x }{\lambda_{1}^q(t)}\right)\,\mathrm{d}x\\
& \quad -\frac{\tilde\rho'(t)}{\eta(t)\lambda_1(t)}\int_{0}^\infty v \psi'\left(\frac{\tilde x }{\lambda_{1}(t)}\right)\varphi \left(\frac{\tilde x }{\lambda_{1}^q(t)}\right)\,\mathrm{d}x  \\ & \quad -\frac{\tilde\rho'(t)}{\eta(t)\lambda_1^q(t)}\int_{0}^\infty v \psi\left(\frac{\tilde x }{\lambda_{1}(t)}\right)\varphi' \left(\frac{\tilde x }{\lambda_{1}^q(t)}\right)\,\mathrm{d}x\\
& \eqqcolon  \Xi_{1,1}[v](t)+\Xi_{1,2}[v](t)+\Xi_{1,3}[v](t) +\Xi_{1,4}[v](t) +\Xi_{1,5}[v](t).
\end{aligned}
\end{align*}
We estimate the terms in this decomposition.

First, by construction, the function $\psi'\left(\frac{\cdot-\tilde{\rho}}{\lambda_1(t) }\right)$ is supported on the set 
\begin{align*}
\big[\tilde\rho(t)-\lambda_{1}(t), \tilde\rho(t)+\lambda_{1}(t)\big]
=\Big[At^{n}-\frac{4Bt^{b}}{\log t}, At^{n}+\frac{4Bt^{b}}{\log t}\Big].
\end{align*}

As a consequence, since $n>b$, for large $t$ this set does not contain the
origin. The same holds for the terms involving $\varphi'$, whose integrand
is supported in $[\tilde\rho(t)-\lambda_{1}^{q}(t),\tilde\rho(t)+\lambda_{1}^{q}(t)]$,
provided $bq<n$. This holds by choosing
$q_{0}<(n-b)/b$. Therefore we can perform integration by
parts without boundary contributions.

Using these support properties and the smoothness of $\varphi=\psi'$, integration by parts gives
\begin{align*}
\begin{aligned}
\Xi_{1,1}[v](t)&=-\frac{1}{\eta(t)}\int_{0}^\infty\partial_x\left(\partial_{xx} v+\frac{1}{2}v^2 \right)\psi\left(\frac{\tilde x }{\lambda_{1}(t)}\right)\varphi\left(\frac{\tilde x }{\lambda^q_{1}(t)}\right)\,\mathrm{d}x\\
&=\frac{1}{\eta(t)\lambda_{1}(t)}\int_{0}^\infty\partial_{xx} v\psi'\left(\frac{\tilde x }{\lambda_{1}(t)}\right)\varphi\left(\frac{ \tilde x }{\lambda^q_{1}(t)}\right)\,\mathrm{d}x \\ & \quad +\frac{1}{\eta(t)\lambda^q_{1}(t)}\int_{0}^\infty \partial_{xx} v\psi\left(\frac{\tilde x }{\lambda_{1}(t)}\right)\varphi'\left(\frac{\tilde x }{\lambda^q_{1}(t)}\right)\,\mathrm{d}x \\
&\quad +\frac{1}{2\eta(t)\lambda_{1}(t)}\int_{0}^\infty  v^2\psi'\left(\frac{\tilde x }{\lambda_{1}(t)}\right)\varphi\left(\frac{\tilde x }{\lambda_{1}^q(t)}\right)\,\mathrm{d}x \\ & \quad +\frac{1}{2\eta(t)\lambda^q_{1}(t)}\int_{0}^\infty  v^2\psi\left(\frac{\tilde x }{\lambda_{1}(t)}\right)\varphi'\left(\frac{ \tilde x }{\lambda^q_{1}(t)}\right)\,\mathrm{d}x\\
& \eqqcolon  \Xi_{1,1,1}[v](t)+\Xi_{1,1,2}[v](t)+\Xi_{1,1,3}[v](t)+\Xi_{1,1,4}[v](t).
\end{aligned}
\end{align*}
We estimate the first term as follows:
\begin{align*}
\begin{aligned}
   \Xi_{1,1,1}[v] &=\frac{1}{\eta(t)\lambda_{1}(t)}\int_{0}^\infty\partial_{xx} v\psi'\left(\frac{\tilde x }{\lambda_{1}(t)}\right)\varphi\left(\frac{ \tilde x }{\lambda^q_{1}(t)}\right)\,\mathrm{d}x \\
   & =-\frac{1}{\eta(t)\lambda_{1}^2(t)}\int_{0}^\infty\partial_x v\psi''\left(\frac{\tilde x }{\lambda_{1}(t)}\right)\varphi\left(\frac{ \tilde x }{\lambda^q_{1}(t)}\right)\,\mathrm{d}x\\
   &\quad - \frac{1}{\eta(t)\lambda_{1}^{q+1}(t)}\int_{0}^\infty\partial_x v\psi'\left(\frac{\tilde x }{\lambda_{1}(t)}\right)\varphi'\left(\frac{ \tilde x }{\lambda^q_{1}(t)}\right)\,\mathrm{d}x\\
&=\frac{1}{\eta(t)\lambda_{1}^3(t)}\int_{0}^\infty v\psi'''\left(\frac{\tilde x }{\lambda_{1}(t)}\right)\varphi\left(\frac{ \tilde x }{\lambda^q_{1}(t)}\right)\,\mathrm{d}x\\
   &\quad + \frac{2}{\eta(t)\lambda_{1}^{q+2}(t)}\int_{0}^\infty v\psi''\left(\frac{\tilde x }{\lambda_{1}(t)}\right)\varphi'\left(\frac{ \tilde x }{\lambda^q_{1}(t)}\right)\,\mathrm{d}x\\
   &\quad + \frac{1}{\eta(t)\lambda_{1}^{2q+1}(t)}\int_{0}^\infty v\psi'\left(\frac{\tilde x }{\lambda_{1}(t)}\right)\varphi''\left(\frac{ \tilde x }{\lambda^q_{1}(t)}\right)\,\mathrm{d}x.
    \end{aligned}
\end{align*}

Using the Cauchy--Schwarz inequality and $q>1$, we obtain
\begin{align*}
\begin{aligned}
    |\Xi_{1,1,1}[v](t)|&\lesssim \frac{1}{\eta(t)\lambda_1^{5/2}(t)}+\frac{1}{\eta(t)\lambda_1^{3/2+q}(t)}+\frac{1}{\eta(t)\lambda_1^{1/2+2q}(t)} \\
    &\lesssim \frac{1}{\eta(t)\lambda_1^{5/2}(t)} \lesssim \frac{1}{t^{m+\frac52b} \log^{-\frac12}t},
\end{aligned}    
\end{align*}
which is integrable in time since $m+\frac52b >1$.

As in \cite[proof of Lemma~5.2]{CavalcanteMunoz2023}, two integrations by parts and the Cauchy--Schwarz inequality give
\begin{align*}
|\Xi_{1,1,2}[v](t)|&\lesssim \frac{1}{\eta(t)\lambda_1^{q+3/2}(t)}+\frac{1}{\eta(t)\lambda_1^{2q+1/2}(t)}+\frac{1}{\eta(t)\lambda_1^{5q/2}(t)},\\
|\Xi_{1,1,4}[v](t)|&\lesssim \frac{1}{\eta(t)\lambda_1^q(t)}.
\end{align*}
Both bounds are integrable on $[T_0,\infty)$ because $m=1-b$ and $q>1$. The term $\Xi_{1,1,3}$ is non-negative and is retained on the left-hand side of the virial inequality.

The translation terms $\Xi_{1,4}(t)$ and $\Xi_{1,5}(t)$ are bounded by
\begin{align*}
|\Xi_{1,4}[v](t)| \lesssim \frac{|\tilde{\rho}'(t)|}{\eta(t)\lambda_1^{1/2}(t)}, \quad |\Xi_{1,5}[v](t)| \lesssim\frac{|\tilde{\rho}'(t)|}{\eta(t)\lambda_1^{q/2}(t)}.
\end{align*} 
Since $q>1$, it is enough to consider the first case. We have
\begin{align*}
|\Xi_{1,4}[v](t)| \lesssim \frac{1}{t^{2-n-b/2} \log^{3/2} t}.
\end{align*}
From \cref{cc1_new}, we conclude that this term is integrable on $[T_0,\infty)$.

To bound  $\Xi_{2}$, we use that $v(t,\cdot)\in L^2(0,\infty)$, uniformly in time, and  that $\lambda_1(t)= 4B t^b /\log t $ and $\eta(t)=t^{m}\ln ^{2} t$: 
\begin{align*}
\begin{aligned}
|\Xi_{2}[v](t)|&=\left|\frac{\eta'(t)}{\eta^{2}(t)}\int_{0}^\infty v\psi\left(\frac{\tilde x }{\lambda_{1}(t)}\right)\varphi\left(\frac{\tilde x }{\lambda^q_{1}(t)}\right)\,\mathrm{d}x\right| \lesssim \frac{(\lambda_1(t))^{q/2}}{t \eta(t)} \lesssim\frac{1}{t^{1+m -\frac{1}{2}bq} \ln^{2 +\frac{q}2} t}.
\end{aligned}
\end{align*}

The Cauchy--Schwarz inequality also gives
\begin{align*}
|\Xi_{1,2}[v](t)|+|\Xi_{1,3}[v](t)|\lesssim \frac{\lambda_1^{1/2}(t)+\lambda_1^{q/2}(t)}{t\eta(t)}\lesssim\frac{\lambda_1^{q/2}(t)}{t\eta(t)}.
\end{align*}
This bound and the bound for $\Xi_2$ are integrable by \cref{cc2_new}. We conclude that
\begin{align*}
\begin{aligned}
\frac{1}{\lambda_{1}(t)\eta(t)}\int_{0}^\infty v^2 \psi'\left(\frac{\tilde x }{\lambda_{1}(t)}\right)\varphi\left(\frac{\tilde x }{\lambda_{1}^q(t)}\right)\,\mathrm{d}x &\le C_0 \frac{\mathrm{d}}{\mathrm{d}t}\Xi[v](t) + \Xi_{\text{int}}(t),
\end{aligned}
\end{align*}
where $\Xi_{\text{int}}\in L^1(T_0,\infty)$ is non-negative and $C_0>0$. The implicit constants depend on the fixed parameters, the cutoffs, and $\|(u_0,v_0,w_0)\|_{L^2(\mathcal Y)}$.

{Since $\psi'\ge c_{0}>0$ on $[-\tfrac12,\tfrac12]$ and
$\varphi\big(\tilde x/\lambda_{1}^{q}(t)\big)\ge c_{0}$ on the support of
$\psi'\big(\tilde x/\lambda_{1}(t)\big)$ for $t$ large (because
$\lambda_{1}^{1-q}(t)\to 0$), and since
$\lambda_{1}(t)/2=2Bt^{b}/\log t\ge Bt^{b}/\log t$, we obtain
\begin{align*}
\int_{T_0}^{\infty}\frac{1}{t\log t}
\int_{At^{n}-\frac{Bt^{b}}{\log t}}^{At^{n}+\frac{Bt^{b}}{\log t}}
v^{2}(t,x)\,\mathrm{d}x\,\mathrm{d}t\lesssim 1 .
\end{align*}
Consequently, since $1/(t\log t)$ is not integrable in time, and recalling
\cref{rho_new},
\begin{align*}
\liminf_{t\to+\infty}\int_{J(t)} v^{2}(t,x)\,\mathrm{d}x=0 .
\end{align*}
The same argument on the other positive edge gives
\begin{align*}
\int_{T_0}^{\infty}\frac{1}{t\log t}\int_{J(t)} w^{2}(t,x)\,\mathrm{d}x\,\mathrm{d}t\lesssim1.
\end{align*}}

{\uline{Step 2:} \emph{Analysis on the negative edge.} On the first edge the equation reads $\partial_t u+\partial_{xxx} u+\gamma_1u\,\partial_x u=0$, and the only
new feature is the coefficient $\gamma_1$, whose sign is not prescribed. We therefore consider the
functional
\begin{align}\label{captal-xi_left}
\Theta[u](t) \coloneqq \frac{\sigma}{\eta(t)}\int_{-\infty}^{0} u(t,x)\,
\psi\left(\frac{x+\tilde\rho(t)}{\lambda_{1}(t)}\right)
\varphi\left(\frac{x+\tilde\rho(t)}{\lambda_{1}^q(t)}\right)\,\mathrm{d}x,
\qquad \sigma\coloneqq\operatorname{sgn}(\gamma_1),
\end{align}
with the same $\lambda_1$, $\eta$, $\tilde\rho$, and $q$ as in \cref{captal-xi}. The weights are centred at
$-\tilde\rho(t)=-At^{n}$, on the negative edge. Setting $\tilde x\coloneqq x+\tilde\rho(t)$, the computation of
$\frac{\mathrm{d}}{\mathrm{d}t}\Theta[u](t)$ has the same structure as in Step~1. The translation terms change sign, which does not affect the estimates of their absolute values. We discuss the nonlinear contribution and the absence of boundary terms.

First, the nonlinear contribution is
\begin{align*}
\begin{aligned}
-\frac{\sigma\gamma_1}{\eta(t)}&\int_{-\infty}^{0}u\,\partial_x u\,
\psi\left(\tfrac{\tilde x}{\lambda_1}\right)\varphi\left(\tfrac{\tilde x}{\lambda_1^q}\right)\,\mathrm{d}x
=\\
&\frac{|\gamma_1|}{2\eta(t)\lambda_1(t)}\int_{-\infty}^{0}u^2\,
\psi'\left(\tfrac{\tilde x}{\lambda_1}\right)\varphi\left(\tfrac{\tilde x}{\lambda_1^q}\right)\,\mathrm{d}x
+\frac{|\gamma_1|}{2\eta(t)\lambda_1^q(t)}\int_{-\infty}^{0}u^2\,
\psi\left(\tfrac{\tilde x}{\lambda_1}\right)\varphi'\left(\tfrac{\tilde x}{\lambda_1^q}\right)\,\mathrm{d}x,
\end{aligned}
\end{align*}
since $\sigma\gamma_1=|\gamma_1|$. Hence the factor $\sigma$ in \cref{captal-xi_left} makes the principal
term coercive regardless of the sign of $\gamma_1$, and it plays the role of $\Xi_{1,1,3}$ in Step~1
(with $\frac12$ replaced by $\frac{|\gamma_1|}{2}$). Positivity of this coefficient uses $\gamma_1\neq0$. All the remaining
terms satisfy the estimates of Step~1, with constants depending on $|\gamma_1|$.

Second, the integrations by parts produce no boundary contributions at the vertex. Since
$\varphi$ is supported in $[-1,1]$, the integrand in \cref{captal-xi_left} and in all the terms of
$\frac{\mathrm{d}}{\mathrm{d}t}\Theta[u](t)$ is supported in
\begin{align*}
\Big[-At^{n}-\Big(\tfrac{4Bt^{b}}{\log t}\Big)^{q},\,-At^{n}+\Big(\tfrac{4Bt^{b}}{\log t}\Big)^{q}\Big],
\end{align*}
and this interval does not contain the origin for $t$ large provided $bq<n$, which holds after choosing
$q_0<(n-b)/b$. In particular the traces $u(t,0)$ and $\partial_x u(t,0)$ do not appear,
so no cancellation between traces on different edges is needed, unlike in \cref{sec:weak}.

Repeating the argument of Step~1 we obtain
\begin{align*}
\int_{T_0}^{\infty}\frac{1}{t\log t}\int_{I(t)} u^{2}(t,x)\,\mathrm{d}x\,\mathrm{d}t\lesssim 1 .
\end{align*}

\uline{Step 3:} \emph{Conclusion.} Adding the three edge estimates and using the definition of $\Omega(t)$, we obtain
\begin{align*}
\int_{T_0}^{\infty}\frac{1}{t\log t}\int_{\Omega(t)}\mathcal U_e^{2}(t,x)\,\mathrm{d}x\,\mathrm{d}t\lesssim 1 .
\end{align*}
Since $\int_{T_0}^{\infty}\frac{\mathrm{d}t}{t\log t}=+\infty$, this forces
$\liminf_{t\to+\infty}\int_{\Omega(t)}\mathcal U_e^{2}(t,x)\,\mathrm{d}x=0$, which concludes the proof of \cref{theo_1}.}

\end{proof}

\section{Analysis of the weakly dominated region}\label{sec:weak}

We prove \cref{teo_2} by estimating solutions in the region $\Omega_0(t)$ near the vertex. Unlike the analysis in \cref{sec:sublinear}, the weights used here reach the vertex $x=0$, so we must cancel the boundary trace terms.
We choose the weight on the negative edge according
to the sign of the nonlinear coefficient $\gamma_1$, and choose
positive constants $A_1,A_2,A_3$ so that the boundary traces cancel.
\begin{proof}
We adapt \cite[Proposition~4.1 and Lemma~4.2]{CavalcanteMunoz2023}, keeping the boundary terms until the three edge identities are added. For $L^2$ solutions, we use the summed identity in the weak formulation. Its three spatial weights vanish together with their first derivatives at the vertex. By \cref{eq:positive-multipliers}, their second derivatives satisfy the adjoint transmission condition, so all vertex terms cancel. The remaining terms contain only the solution and its square against smooth decaying weights, which are integrable by the $L^2$ bound. The edgewise integrations by parts in this proof display this cancellation for regular solutions.

Throughout this section, we assume $\gamma_1\neq0$ and set $\sigma\coloneqq\operatorname{sgn}(\gamma_1)\in\{-1,1\}$. We also assume that there exist constants $A_1,A_2,A_3>0$ such that
\begin{align}\label{eq:boundary-cancellation}
    \sigma A_1
    +\frac{A_2}{\alpha_2}
    +\frac{A_3}{\alpha_3}
    =0.
\end{align}

The condition \cref{eq:boundary-cancellation} is homogeneous in the constants
$A_j$. Thus, only their relative sizes are relevant. In particular,
we may normalize $A_2=1$ or $A_3=1$.
 We use
\begin{align*}
    \psi_+(s)\coloneqq\tanh(s^2),
    \qquad
    \psi_-(s)\coloneqq-\sigma\tanh(s^2),
    \qquad
    \varphi(s)\coloneqq\operatorname{sech} s.
\end{align*}
These functions are smooth and satisfy
\begin{align*}
    \psi_+(0)=\psi_-(0)=0,
    \qquad
    \psi_+'(0)=\psi_-'(0)=0.
\end{align*}
Moreover,
\begin{align*}
    \psi_+'(s)
    =
    2s\operatorname{sech}^2(s^2)>0,
    \qquad s>0.
\end{align*}
On the negative half-line,
\begin{align*}
    \psi_-'(s)
    =
    -2\sigma s\operatorname{sech}^2(s^2),
\end{align*}
so that
\begin{align}\label{eq:gamma-psi-positive}
    \gamma_1\psi_-'(s)
    =
    -2|\gamma_1|s\operatorname{sech}^2(s^2)>0,
    \qquad s<0.
\end{align}
Thus, the nonlinear contribution on the negative edge has the same
sign as the corresponding contributions on the two positive edges,
regardless of the sign of $\gamma_1$.

We define the functions

\begin{align}\label{eq-lambda-eta}
\lambda_{1}(t)=\frac{t^{b}}{\ln t},\qquad \eta(t)=t^{m}\ln^{2}t,
\end{align}
and assume the conditions
\begin{align}
    &m,b>0,\qquad b\leq 1-m, \label{eq:cond1}\\
    &\frac25(1-m)<b\leq\frac12(m+1), \label{eq:cond2}\\
    &\frac1q(1-m)<b\leq\frac{2m}{q}, \label{eq:cond3}\\
    &b\leq\frac12(1-m). \label{eq:cond4}
\end{align}

{
Given $0<\varepsilon_1<\varepsilon_2$ as in the statement of \cref{teo_2}, we fix
\begin{align}\label{eq:params}
  q=\tfrac52,\qquad \delta_0\coloneqq\tfrac12(\varepsilon_1+\varepsilon_2),\qquad
  m=\tfrac{5}{13}+\delta_0,\qquad b=\tfrac{4}{13}-\delta_0,
\end{align}
so that $\varepsilon_1<\delta_0<\varepsilon_2$. The restriction $\varepsilon_2<4/39$ ensures that conditions
\cref{eq:cond1,eq:cond2,eq:cond3,eq:cond4} hold, and moreover $m+\tfrac52 b=\tfrac{15}{13}-\tfrac32\delta_0>1$
and $2b+m=1-\delta_0$.}

\medskip
\noindent

\uline{Step 1:} \emph{Virial functionals.} For the two positive edges, we define
\begin{align*}
    \Xi[v](t)
    \coloneqq
    \frac{1}{\eta(t)}
    \int_0^\infty
    v(t,x)
    \psi_+\left(\frac{x}{\lambda_1(t)}\right)
    \varphi\left(\frac{x}{\lambda_1^q(t)}\right)\,\mathrm{d}x
\end{align*}
and
\begin{align*}
    \Xi[w](t)
    \coloneqq
    \frac{1}{\eta(t)}
    \int_0^\infty
    w(t,x)
    \psi_+\left(\frac{x}{\lambda_1(t)}\right)
    \varphi\left(\frac{x}{\lambda_1^q(t)}\right)\,\mathrm{d}x.
\end{align*}
On the negative edge, we define
\begin{align}\label{eq:Xi-u-new}
    \Xi[u](t)
    \coloneqq
    \frac{1}{\eta(t)}
    \int_{-\infty}^0
    u(t,x)
    \psi_-\left(\frac{x}{\lambda_1(t)}\right)
    \varphi\left(\frac{x}{\lambda_1^q(t)}\right)\,\mathrm{d}x.
\end{align}
Finally, we consider the combined functional
\begin{align}\label{eq:Xi-total-new}
    \Xi_{\mathrm{tot}}(t)
    \coloneqq
    A_1\Xi[u](t)
    +A_2\Xi[v](t)
    +A_3\Xi[w](t).
\end{align}

Since the functions $\psi_\pm$ and $\varphi$ are bounded, the
Cauchy--Schwarz inequality gives
\begin{align*}
    |\Xi[v](t)|
    \lesssim
    \frac{\lambda_1^{q/2}(t)}{\eta(t)}
    \|v(t,\cdot)\|_{L^2(0,\infty)},
\end{align*}
and the same estimate holds for $\Xi[w]$ and $\Xi[u]$.
Using \cref{eq-lambda-eta} and the uniform $L^2$ bound in
\cref{cor_1}, we obtain
\begin{align}\label{eq:Xi-bounded}
    \sup_{t\gg1}
    \left(
    |\Xi[u](t)|+
    |\Xi[v](t)|+
    |\Xi[w](t)|
    \right)<\infty.
\end{align}

\medskip
\noindent
\uline{Step 2:} \emph{The positive edges.} We first consider $\Xi[v]$. 

We estimate the time derivative of $\Xi[v](t)$, retaining the boundary contribution in \cref{first_edge_positive} until the three edge estimates are added. We have
\begin{align*}
\begin{aligned}
 \frac{\mathrm{d}}{\mathrm{d}t}\Xi[v](t)&= \frac{1}{\eta(t)}\int_{0}^\infty\partial_t\left(v\psi_{+}\left(\frac{x }{\lambda_{1}(t)}\right)\varphi\left(\frac{x }{\lambda_{1}^q(t)}\right)\right)\,\mathrm{d}x \\ &\qquad -\frac{\eta'(t)}{\eta^{2}(t)}\int_{0}^\infty v\psi_{+}\left(\frac{x}{\lambda_{1}(t)}\right)\varphi\left(\frac{x }{\lambda_{1}^q(t)}\right)\,\mathrm{d}x\\
& \eqqcolon  \Xi_{1}[v](t)+ \Xi_{2}[v](t).
\end{aligned}
\end{align*}

Using the uniform $L^2$ bound in \cref{cor_1} and the definitions in \cref{eq-lambda-eta}, we obtain
\begin{align*}
\begin{aligned}
|\Xi_{2}[v](t)|&\leq\left|\frac{\eta'(t)}{\eta^{2}(t)}\int_{0}^\infty v\psi_{+}\left(\frac{x }{\lambda_{1}(t)}\right)\varphi\left(\frac{x }{\lambda^q_{1}(t)}\right)\,\mathrm{d}x\right| \lesssim \frac{(\lambda_1(t))^{q/2}}{t \eta(t)} =\frac{1}{t^{1+m -\frac{1}{2}bq} \ln^{2 +\frac{q}2} t}.
\end{aligned}
\end{align*}
By \cref{eq:cond3}, this bound is integrable in time. Thus,  $\Xi_{2}\in L^1(\{t\gg1\})$.

We expand $\Xi_1[v]$ to identify its boundary contribution:
\begin{align*}
\begin{aligned}
\Xi_{1}[v](t)&=\frac{1}{\eta(t)}\int_{0}^\infty\partial_t v\psi_{+}\left(\frac{x }{\lambda_{1}(t)}\right)\varphi\left(\frac{x }{\lambda_{1}^q(t)}\right)\,\mathrm{d}x \\
&\quad -\frac{\lambda_{1}'(t)}{\lambda_{1}(t)\eta(t)}\int_{0}^\infty v \psi_{+}'\left(\frac{x }{\lambda_{1}(t)}\right)\left(\frac{x }{\lambda_{1}(t)}\right)\varphi\left(\frac{x }{\lambda_{1}^q(t)}\right)\,\mathrm{d}x\\
&\quad -q\frac{\lambda_{1}'(t)}{\lambda_{1}(t)\eta(t)}\int_{0}^\infty v \psi_{+}\left(\frac{x }{\lambda_{1}(t)}\right)\left(\frac{x }{\lambda_{1}^q(t)}\right)\varphi'\left(\frac{x }{\lambda_{1}^q(t)}\right)\,\mathrm{d}x\\
& \eqqcolon  \Xi_{1,1}[v](t)+\Xi_{1,2}[v](t)+\Xi_{1,3}[v](t) .
\end{aligned}
\end{align*}
Using $\psi_{+}(0)=0$ and the rapid decay of $\varphi$, integration by parts gives
\begin{align}\label{separation_Xi1_gKdV3}
\begin{aligned}
\Xi_{1,1}[v](t)&=-\frac{1}{\eta(t)}\int_{0}^\infty\partial_x\left(\partial_{xx} v+\frac{v^2}{2} \right)\psi_{+}\left(\frac{x }{\lambda_{1}(t)}\right)\varphi\left(\frac{x }{\lambda^q_{1}(t)}\right)\,\mathrm{d}x\\
&=\frac{1}{\eta(t)\lambda_{1}(t)}\int_{0}^\infty\partial_{xx} v\psi_{+}'\left(\frac{x }{\lambda_{1}(t)}\right)\varphi\left(\frac{x }{\lambda^q_{1}(t)}\right)\,\mathrm{d}x \\&\quad+\frac{1}{\eta(t)\lambda^q_{1}(t)}\int_{0}^\infty \partial_{xx} v\psi_{+}\left(\frac{x }{\lambda_{1}(t)}\right)\varphi'\left(\frac{x }{\lambda^q_{1}(t)}\right)\,\mathrm{d}x \\
&\quad +\frac{1}{2\eta(t)\lambda_{1}(t)}\int_{0}^\infty  v^2\psi_{+}'\left(\frac{x }{\lambda_{1}(t)}\right)\varphi\left(\frac{x }{\lambda_{1}^q(t)}\right)\,\mathrm{d}x\\&\quad +\frac{1}{2\eta(t)\lambda^q_{1}(t)}\int_{0}^\infty  v^2\psi_{+}\left(\frac{x }{\lambda_{1}(t)}\right)\varphi'\left(\frac{x }{\lambda^q_{1}(t)}\right)\,\mathrm{d}x\\
& \eqqcolon  \Xi_{1,1,1}[v](t)+\Xi_{1,1,2}[v](t)+\Xi_{1,1,3}[v](t)+\Xi_{1,1,4}[v](t).
\end{aligned}
\end{align}
We bound the time integral of $\Xi_{1,1,3}[v]$. 

Using $\psi_{+}'(0)=0$ and the rapid decay of $\varphi$, further integration by parts gives
\begin{align*}
\begin{aligned}
\Xi_{1,1,1}[v](t) & =-\frac{1}{\eta(t)\lambda_{1}(t)}\int_{0}^\infty\partial_x v \,\partial_x \left( \psi_{+}'\left(\frac{x }{\lambda_{1}(t)}\right)\varphi\left(\frac{x }{\lambda^q_{1}(t)}\right) \right)\,\mathrm{d}x\\
& =\frac{1}{\eta(t)\lambda_{1}(t)}\int_{0}^\infty v \,\partial_{xx} \left( \psi_{+}'\left(\frac{x }{\lambda_{1}(t)}\right)\varphi\left(\frac{x }{\lambda^q_{1}(t)}\right) \right)\,\mathrm{d}x+\mathcal T(v)(t),
\end{aligned}
\end{align*}
where the trace term is given by
\begin{align}\label{trace_right_1}
\begin{aligned}
\mathcal T (v)(t) &\coloneqq \frac{1}{\eta(t)\lambda_1(t)}v(t,0)\partial_x \left( \psi_{+}'\left(\frac{x }{\lambda_{1}(t)}\right)\varphi\left(\frac{x }{\lambda^q_{1}(t)} \right)\right)\bigg|_{x=0}\\&= \frac{v(t,0)}{\eta(t)} \left( \frac{\psi_{+}''(0)\varphi(0)}{\lambda_1^2(t)} + \frac{\psi_{+}'(0)\varphi'(0)}{\lambda_1^{q+1}(t)} \right)\\
&= \frac{2}{\lambda_1^2(t)\eta(t)} v(t,0),
\end{aligned}
\end{align}
where we have used that  $\psi_+''(0)=2$ and $\psi_+'(0)=0$.
Using \cref{separation_Xi1_gKdV3,trace_right_1}, we obtain 
\begin{align*}
\begin{aligned}
\Xi_{1,1,1}[v](t) & =\frac{1}{\eta(t)\lambda_{1}^{3}(t)}\int_{0}^\infty v\psi_{+}'''
\left(\frac{x }{\lambda_{1}(t)}\right)\varphi\left(\frac{x }{\lambda^q_{1}(t)}\right)\,\mathrm{d}x \\
& \quad +\frac{2}{\eta(t)\lambda_{1}^{2+q}(t)}\int_{0}^\infty v\psi_{+}''\left(\frac{x }{\lambda_{1}(t)}\right)\varphi'\left(\frac{x }{\lambda^q_{1}(t)}\right)\,\mathrm{d}x\\
&\quad +\frac{1}{\eta(t)\lambda_{1}^{1+2q}(t)}\int_{0}^\infty v\psi_{+}'\left(\frac{x }{\lambda_{1}(t)}\right)\varphi''\left(\frac{x }{\lambda^q_{1}(t)}\right)\,\mathrm{d}x\\
&\quad+ \frac{2}{\lambda_1^2(t)\eta(t)} v(t,0).
\end{aligned}
\end{align*}
First, we bound each term of $\Xi_{1,1,1}[v]$ using the Cauchy--Schwarz inequality and the fact that $q>1$: 
\begin{align*}
\begin{aligned}
    |\Xi_{1,1,1}[v](t)-\frac{2}{\lambda_1^2(t)\eta(t)} v(t,0)|&\lesssim \frac{1}{\eta(t)\lambda_1^{5/2}(t)}+\frac{1}{\eta(t)\lambda_1^{3/2+q}(t)}+\frac{1}{\eta(t)\lambda_1^{1/2+2q}(t)} \\
    &\lesssim \frac{1}{\eta(t)\lambda_1^{5/2}(t)} \lesssim \frac{1}{t^{m+\frac52b} \log^{-\frac12}t}.
\end{aligned}    
\end{align*}
This bound is integrable in time since $m+\frac52b >1$.

We estimate $\Xi_{1,1,2}[v](t)$ and $\Xi_{1,1,4}[v](t)$ as in \cite[proof of Lemma~4.2]{CavalcanteMunoz2023}:
\begin{align*}
  \begin{aligned}
    |\Xi_{1,1,2}[v](t)|     &\lesssim \frac{1}{t^{m+(\frac32+q)b} \log^{\frac12 -q}t}.
  \end{aligned}
\end{align*}

\begin{align*}
\begin{aligned}
|\Xi_{1,1,4}[v](t)| &\lesssim \frac{1}{2t^{m+bq} \log^{2-q} t}.
\end{aligned}
\end{align*}
Since $q>\frac{1-m}{b}$, $\Xi_{1,1,4}\in L^1(\{t\gg1\})$.

By the Cauchy--Schwarz inequality and the uniform $L^2$ bound, we obtain
\begin{align*}
\left|\Xi_{1,2}[v](t)\right|
&\lesssim \frac{|\lambda_1'(t)|}{\eta(t)\lambda_1^{1/2}(t)}
\lesssim \frac{1}{t^{1+m-b/2}\ln^{5/2}t}.
\end{align*}
Since $q>1$ and $b\le 2m/q$, we have $b/2<m$, so $\Xi_{1,2}[v]\in L^1(\{t\gg1\})$.

We consider the term $\Xi_{1,3}[v](t)$. Using the properties of $\varphi$, the uniform $L^2$ bound for $v$, and Young's inequality with $\theta(t)=t^m$, we obtain
\begin{align*}
\begin{aligned}
\left| \Xi_{1,3}[v](t) \right|   &\lesssim  \left|\frac{\lambda_{1}'(t)\theta(t)}{\lambda_{1}(t)\eta(t)}\right|   + \left|\frac{\lambda_{1}'(t)\lambda_1^{q}(t)}{\lambda_1(t)\eta(t)\theta(t)}\right|.
\end{aligned}
\end{align*}
Hence,
\begin{align*}
\begin{aligned}
\left| \Xi_{1,3}[v](t) \right|   &\lesssim \frac{1}{t\ln^2t}+\frac{1}{t^{1+2m -bq} \ln^{2+q} t},
\end{aligned}
\end{align*}
which is integrable in time if $b\leq \frac{2m}{q}$. Thus, $\Xi_{1,3}[v]\in L^1(\{t\gg1\})$. We conclude that
\begin{align}\label{first_edge_positive}
\begin{aligned}
\frac{1}{2\lambda_{1}(t)\eta(t)}\int_{0}^\infty v^2 \psi_{+}'\left(\frac{x }{\lambda_{1}(t)}\right)\varphi\left(\frac{x }{\lambda_{1}^q(t)}\right)\,\mathrm{d}x +  \frac{2}{\lambda_1^2(t)\eta(t)} v(t,0)&\le \frac{\mathrm{d}}{\mathrm{d}t}\Xi[v](t) + \Xi_{\text{int}}[v](t),
\end{aligned}
\end{align}
where $\Xi_{\text{int}}[v]$ is an integrable remainder for sufficiently large $t$.

By performing the same computations, we obtain a similar estimate for the third edge:

\begin{align}\label{second_edge_positive}
\begin{aligned}
\frac{1}{2\lambda_{1}(t)\eta(t)}\int_{0}^\infty w^2 \psi_{+}'\left(\frac{x }{\lambda_{1}(t)}\right)\varphi\left(\frac{x }{\lambda_{1}^q(t)}\right)\,\mathrm{d}x +  \frac{2}{\lambda_1^2(t)\eta(t)} w(t,0)&\le \frac{\mathrm{d}}{\mathrm{d}t}\Xi[w](t) + \Xi_{\text{int}}[w](t).
\end{aligned}
\end{align}
\uline{Step 3:} \emph{The negative edge.} We consider the functional $\Xi[u]$ defined in
\cref{eq:Xi-u-new}. The equation on the first edge is
\begin{align*}
    \partial_t u+\partial_{xxx} u+\gamma_1u\,\partial_x u=0.
\end{align*}

The computation is the same as on the positive edges, except that the coefficient of the nonlinearity is
$\gamma_1$ and that the weight is $\psi_-$.

The nonlinear contribution is
\begin{align*}
    &-\frac{\gamma_1}{\eta(t)}
    \int_{-\infty}^0
    u\,\partial_x u
    \psi_-\left(\frac{x}{\lambda_1}\right)
    \varphi\left(\frac{x}{\lambda_1^q}\right)\,\mathrm{d}x
    \\
    &\qquad
    =
    \frac{\gamma_1}{2\lambda_1(t)\eta(t)}
    \int_{-\infty}^0
    u^2
    \psi_-'\left(\frac{x}{\lambda_1(t)}\right)
    \varphi\left(\frac{x}{\lambda_1^q(t)}\right)\,\mathrm{d}x
    +\mathcal R_u(t),
\end{align*}
where $\mathcal R_u$ contains the terms in which derivatives fall on
the localized factor $\varphi$.

By \cref{eq:gamma-psi-positive},
\begin{align*}
    \gamma_1
    \psi_-'\left(\frac{x}{\lambda_1(t)}\right)>0,
    \qquad x<0.
\end{align*}
More precisely,
\begin{align*}
    \gamma_1
    \psi_-'\left(\frac{x}{\lambda_1(t)}\right)
    =
    \frac{2|\gamma_1||x|}{\lambda_1(t)}
  \operatorname{sech}^2\left(\frac{x^2}{\lambda_1^2(t)}\right)>0.
\end{align*}
Hence the principal nonlinear contribution on the negative edge is
coercive regardless of the sign of $\gamma_1$.

We next consider the boundary term generated by the dispersive
operator. As in the calculation on the positive edges, integration by parts
gives
\begin{align*}
   \mathcal T_u(t)
    =
    -\frac{\psi_-''(0)}
    {\lambda_1^2(t)\eta(t)}
    u(t,0).
\end{align*}
Since $\psi_-''(0)=-2\sigma$, we obtain
\begin{align*}
   \mathcal T_u(t)
    =
    \frac{2\sigma}
    {\lambda_1^2(t)\eta(t)}
    u(t,0).
\end{align*}

The remaining terms generated by derivatives of $\psi_-$ and $\varphi$ are
estimated exactly as on the positive edges. Since $\psi_-$ is only
a constant multiple of $\psi_+$, all estimates remain unchanged up to
multiplicative constants. Therefore,
\begin{align}\label{eq:virial-u-final}
\begin{aligned}
&\frac{|\gamma_1|}{2\lambda_1(t)\eta(t)}
    \int_{-\infty}^0
    u^2(t,x)
    \left|
    \psi_-'\left(\frac{x}{\lambda_1(t)}\right)
    \right|
    \varphi\left(\frac{x}{\lambda_1^q(t)}\right)\,\mathrm{d}x
    \\
    &\qquad
    +\frac{2\sigma}{\lambda_1^2(t)\eta(t)}
    u(t,0)
    \leq
    \frac{\mathrm{d}}{\mathrm{d}t}\Xi[u](t)
    +\Xi_{\mathrm{int}}[u](t),
\end{aligned}
\end{align}
where $\Xi_{\mathrm{int}}[u]\in L^1(\{t\gg1\})$.

\medskip
\noindent
\uline{Step 4:} \emph{Cancellation of the boundary traces.} We multiply \cref{eq:virial-u-final},
\cref{first_edge_positive}, and
\cref{second_edge_positive} by $A_1$, $A_2$, and $A_3$, respectively,
and add the resulting inequalities.

The boundary contribution is
\begin{align*}
    \mathcal B(t)
    &\coloneqq
    \frac{2}{\lambda_1^2(t)\eta(t)}
    \left[
        \sigma A_1u(t,0)
        +A_2v(t,0)
        +A_3w(t,0)
    \right].
\end{align*}

The transmission conditions at the vertex are
\begin{align*}
    u(t,0)=\alpha_2v(t,0)
    =\alpha_3w(t,0).
\end{align*}
Consequently,
\begin{align*}
    v(t,0)=\frac{u(t,0)}{\alpha_2},
    \qquad
    w(t,0)=\frac{u(t,0)}{\alpha_3},
\end{align*}
and therefore
\begin{align*}
    \mathcal B(t)
    &=
    \frac{2u(t,0)}
    {\lambda_1^2(t)\eta(t)}
    \left(
        \sigma A_1
        +\frac{A_2}{\alpha_2}
        +\frac{A_3}{\alpha_3}
    \right).
\end{align*}
By \cref{eq:boundary-cancellation}, $\mathcal B(t)=0$. We obtain
\begin{align}\label{eq:total-virial}
\begin{aligned}
&\frac{A_1|\gamma_1|}
    {2\lambda_1(t)\eta(t)}
    \int_{-\infty}^0
    u^2(t,x)
    \left|
    \psi_-'\left(\frac{x}{\lambda_1(t)}\right)
    \right|
    \varphi\left(\frac{x}{\lambda_1^q(t)}\right)\,\mathrm{d}x
    \\
    &\quad+
    \frac{A_2}
    {2\lambda_1(t)\eta(t)}
    \int_0^\infty
    v^2(t,x)
    \psi_+'\left(\frac{x}{\lambda_1(t)}\right)
    \varphi\left(\frac{x}{\lambda_1^q(t)}\right)\,\mathrm{d}x
    \\
    &\quad+
    \frac{A_3}
    {2\lambda_1(t)\eta(t)}
    \int_0^\infty
    w^2(t,x)
    \psi_+'\left(\frac{x}{\lambda_1(t)}\right)
    \varphi\left(\frac{x}{\lambda_1^q(t)}\right)\,\mathrm{d}x
    \\
    &\leq
    \frac{\mathrm{d}}{\mathrm{d}t}\Xi_{\mathrm{tot}}(t)
    +\Xi_{\mathrm{int,tot}}(t),
\end{aligned}
\end{align}
where
\begin{align*}
    \Xi_{\mathrm{int,tot}}
    \coloneqq
    A_1\Xi_{\mathrm{int}}[u]
    +A_2\Xi_{\mathrm{int}}[v]
    +A_3\Xi_{\mathrm{int}}[w]
    \in L^1(\{t\gg1\}).
\end{align*}

\medskip

{\medskip\noindent\uline{Step 5:} \emph{Time integration and conclusion.} We fix $t_0>1$ sufficiently large that the estimates in this proof hold. Integrating \cref{eq:total-virial} over $[t_0,T]$, using \cref{eq:Xi-bounded} and the integrability of $\Xi_{\rm int,tot}$, and applying the monotone convergence theorem as $T\to\infty$, we obtain
\begin{align}\label{eq:integrated}
  \int_{t_0}^{\infty}\frac{1}{\lambda_1(t)\eta(t)}
  \bigg[\int_{-\infty}^{0}u^2\,\Big|\psi_-'\Big(\tfrac{x}{\lambda_1}\Big)\Big|\,\varphi\Big(\tfrac{x}{\lambda_1^{q}}\Big)\,\mathrm{d}x
  +\int_{0}^{\infty}\big(v^2+w^2\big)\,\psi_+'\Big(\tfrac{x}{\lambda_1}\Big)\varphi\Big(\tfrac{x}{\lambda_1^{q}}\Big)\,\mathrm{d}x\bigg]\,\mathrm{d}t<\infty .
\end{align}

We derive a lower bound for the weight on $\Omega_{0}(t)$. Let
$x\in J_0(t)=(t^{-\varepsilon_{1}},t^{4/13-\varepsilon_{2}})$ and set
$s\coloneqq x/\lambda_{1}(t)$. Since $\lambda_{1}(t)=t^{4/13-\delta_{0}}/\ln t$ and
$\varepsilon_{2}>\delta_{0}$,
\begin{align*}
0<s&\le \frac{t^{4/13-\varepsilon_{2}}\,\ln t}{t^{4/13-\delta_{0}}}\\
&=t^{-(\varepsilon_{2}-\delta_{0})}\ln t\longrightarrow0,
\qquad t\to+\infty,
\end{align*}
so that for $t$ large we have $s\le\frac12$,
$\operatorname{sech}^{2}(s^{2})\ge\frac12$, and
$\varphi(x/\lambda_{1}^{q})=\operatorname{sech}(s\,\lambda_{1}^{1-q})\ge\frac12$.
Using $\psi_{+}'(s)=2s\operatorname{sech}^{2}(s^{2})\ge s$ for $0<s\le\frac12$ and $x\ge t^{-\varepsilon_{1}}$, we obtain, for all $x\in J_0(t)$ and sufficiently large $t$,
\begin{align}\label{eq:weight-below}
\begin{aligned}
\frac{1}{\lambda_{1}\eta}\,\psi_{+}'\Big(\frac{x}{\lambda_{1}}\Big)
\varphi\Big(\frac{x}{\lambda_{1}^{q}}\Big)
&\ge\frac{x}{2\lambda_{1}^{2}(t)\eta(t)}\\
&=\frac{x}{2t^{2b+m}}
\ge\frac{t^{-1+(\delta_{0}-\varepsilon_{1})}}{2},
\end{aligned}
\end{align}
where we used $2b+m=1-\delta_{0}$ and
$\lambda_{1}^{2}\eta=t^{2b+m}$. Since $|\psi_{-}'|=\psi_{+}'(|\cdot|)$, the
same bound holds for $x\in I_0(t)$. Inserting \cref{eq:weight-below} into
\cref{eq:integrated} yields
\begin{align}\label{eq:final-int}
\int_{t_{0}}^{\infty} t^{-1+(\delta_{0}-\varepsilon_{1})}
\int_{\Omega_{0}(t)} \mathcal U_{e}^{2}(t,x)\,\mathrm{d}x\,\mathrm{d}t<\infty .
\end{align}

Since $\delta_{0}>\varepsilon_{1}$, we have
$\int_{t_{0}}^{\infty}t^{-1+(\delta_{0}-\varepsilon_{1})}\,\mathrm{d}t=+\infty$,
and \cref{eq:final-int} forces
\begin{align*}\liminf_{t\to+\infty}\int_{\Omega_{0}(t)}\mathcal U_{e}^{2}(t,x)\,\mathrm{d}x=0.\end{align*}
This proves the desired decay in the sense of $\liminf$.
}

\end{proof}

\section{General star graphs}
\label{sec:general}
We discuss a formal extension of the arguments in \cref{sec:sublinear,sec:weak} to KdV equations on general star graphs, such as those considered in \cite{CavalcanteMarques2025}. We identify the mass balance and parameter conditions needed for these arguments, assuming an appropriate well-posedness theory.

We consider the KdV equation on a star graph $\mathcal G$ with $k$ left half-lines and $m$ right half-lines: 
\begin{align*}
\begin{cases}
\partial_t u_i+\partial_{xxx} u_i+\gamma_iu_i\partial_x u_i=0,& (t,x)\in (0,T)\times(-\infty,0), i=1,\dots, k,\\
\partial_t v_j+\partial_{xxx} v_j+\gamma_{k+j}v_j\partial_x v_j=0,& (t,x)\in (0,T)\times(0,\infty), j=1,\dots, m,\\
u_i(0,x)=u_{0,i}(x),&  x \in  (-\infty,0), i=1,\dots, k,\\
v_j(0,x)=v_{0,j}(x),&  x \in  (0,\infty), j=1,\dots, m,\\
\end{cases}
\end{align*}
We impose the boundary conditions
\begin{align*}
\begin{cases}
u_1(t,0)=\alpha_i u_i(t,0)=\alpha_{k+j} v_j(t,0), i=2, \dots, k, j=1, \dots, m,\\
\partial_x U(t,0)=B\partial_x V(t,0),\\
\displaystyle\sum^k_{i=1} \alpha_i^{-1}\partial_{xx} u_i(t,0)=\displaystyle\sum^m_{j=1} \alpha^{-1}_{k+j} \partial_{xx} v_j(t,0),
\end{cases}
\end{align*}
where $V(t,x)=(v_1, \dots, v_{m})^T$, $U(t,x)=(u_1, \dots, u_{k})^T$, and $B=\left[b_{ij}\right]_{k\times m}$ is a matrix such that \begin{align}\label{B}\| B \mathbf{x} \|\leq\|\mathbf{x} \|, \quad \text{for all }\mathbf{x} \in \mathbb{R}^m,\end{align} where the norms are Euclidean, $\alpha_1\coloneqq1$, and $\alpha_i,\alpha_{k+j}$ are non-zero constants.

\begin{figure}[htp]
    \centering 
    \begin{tikzpicture}[scale=2]
        \node at (0.8,0.1)[]{\vdots};
        \node at (-0.8,0.1)[]{\vdots};
        
        \foreach \angle in {40,70,290,320} {
            \draw[dashed, black!55] (\angle:2cm) -- (0,0);
            \draw[->,thick,black!55] (0,0) -- (\angle:1.2cm);
        }
        \foreach \angle in {110,140,220,250} {
            \draw[dashed, black!85] (\angle:2cm) -- (0,0);
            \draw[->,thick,black!85] (\angle:1.2cm) -- (\angle:1.1cm);
            \draw[-,thick, black!85] (\angle:1.2cm) -- (0,0);
        }
        \draw[fill] (0,0) circle (1.5pt);
    \end{tikzpicture}
    \caption{A star graph with $m+k$ edges.}
    \label{figure1}
\end{figure}
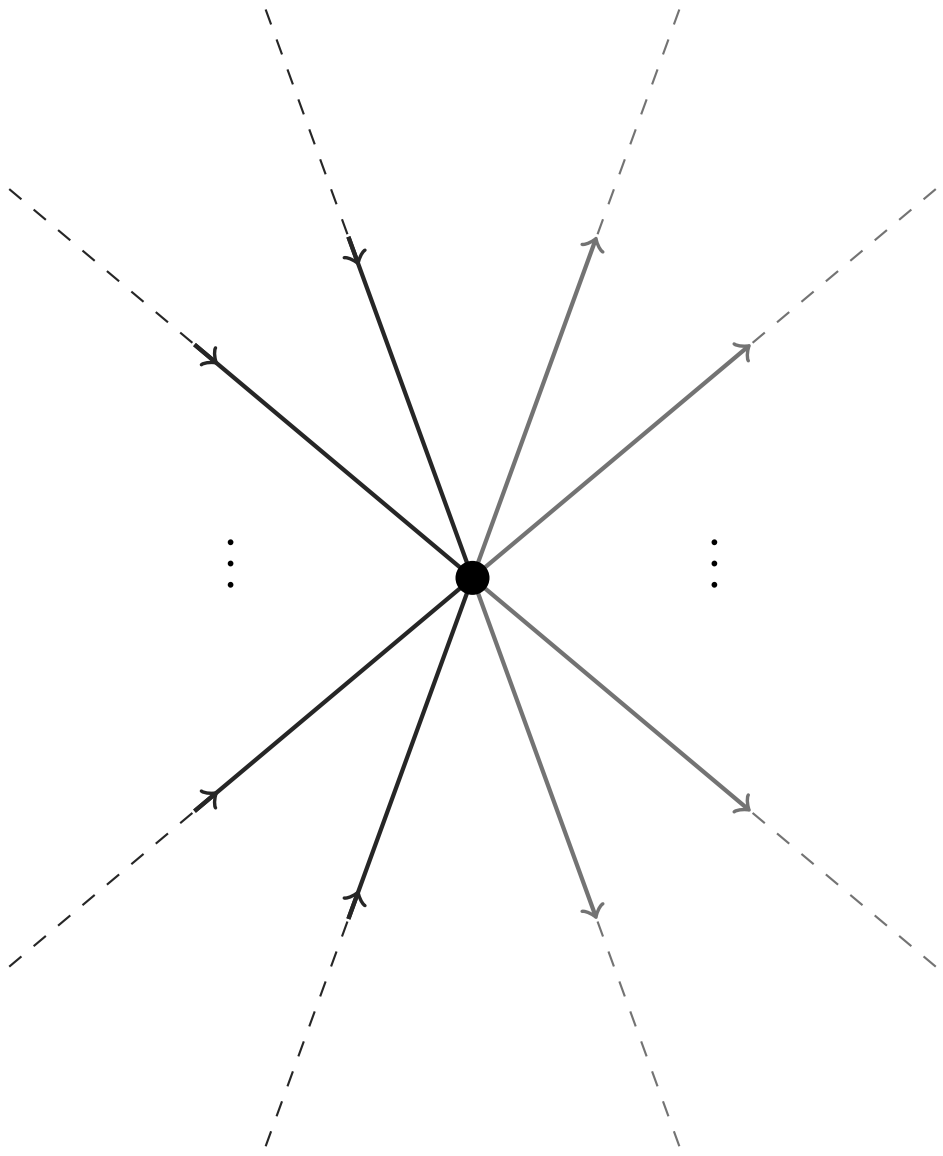

The estimates in \crefrange{sec:sublinear}{sec:weak} require a uniform $L^2$ bound. We therefore compute the mass balance at the vertex.

As in the proof of \cref{zamp1}, integration by parts gives the formal identity
\begin{align*}
\begin{aligned}
\frac12\frac{\mathrm{d}}{\mathrm{d}t} & \left( \sum_{i=1}^k \int_{-\infty}^{0} u_i^2(t,x)\,\mathrm{d}x + \sum_{j=1}^m \int_0^{\infty} v_j^2(t,x)\,\mathrm{d}x \right)\\
&= \underbrace{\left(-\sum_{i=1}^k u_i(t,0)\partial_{xx} u_i(t,0) + \sum_{j=1}^m\partial_{xx} v_j(t,0) v_j(t,0)\right)}_{A} \\
&+ \underbrace{\frac12\left(\sum_{i=1}^k (\partial_x u_i(t,0))^2 - \sum_{j=1}^m(\partial_x v_j(t,0))^2 \right)}_{B}  \\ &\qquad\underbrace{ - \frac13\left( \sum_{i=1}^k\gamma_i u_i^3(t,0) - \sum_{j=1}^m\gamma_{k+j}v_j^3(t,0) \right)}_{C}.
\end{aligned}
\end{align*}

Using the boundary conditions and \cref{B}, we have $A=0$, $B \leq 0$. Furthermore, applying the first boundary condition, we obtain:
\begin{align*}
    C = -\frac{u_1^3(t,0)}{3} \left( \gamma_1 + \sum_{i=2}^k \frac{\gamma_i}{\alpha_i^3} - \sum_{j=1}^m \frac{\gamma_{k+j}}{\alpha_{k+j}^3} \right).
\end{align*}

Thus, if the parameters satisfy
\begin{align}\label{cond_final}
\gamma_1 + \sum_{i=2}^k \frac{\gamma_i}{\alpha_i^3} - \sum_{j=1}^m \frac{\gamma_{k+j}}{\alpha_{k+j}^3} = 0,
\end{align}
the formal identity yields a non-increasing total $L^2$ mass. To extend \cref{theo_1}, we also require $\gamma_e\neq0$ on every edge, a local well-posedness theory with a lifespan controlled by the $L^2$ norm, and a justification of this mass inequality for the resulting solutions. The local theory in \cite[Theorem~2.1]{CavalcanteMarques2025} also requires a certain matrix to be invertible. Subject to these requirements, the edgewise argument in \cref{sec:sublinear} gives decay in the strongly dominated linear region.

To extend the boundary cancellation in \cref{teo_2} to the weakly dominated linear region, we additionally require constants $A_1,\ldots,A_{k+m}>0$ such that
\begin{align}\label{eq:star-boundary-cancellation}
\sum_{i=1}^k\frac{\operatorname{sgn}(\gamma_i)A_i}{\alpha_i}+\sum_{j=1}^m\frac{\operatorname{sgn}(\gamma_{k+j})A_{k+j}}{\alpha_{k+j}}=0.
\end{align}
The signs on both incoming and outgoing edges are needed when the coefficients $\gamma_e$ have arbitrary signs. If global $L^2$ solutions satisfy the mass inequality and $\gamma_e\neq0$ on every edge, \cref{cond_final,eq:star-boundary-cancellation} allow us to apply the argument in \cref{sec:weak} and obtain decay analogous to \cref{teo_2}.

\section*{Acknowledgments}

M.~Cavalcante acknowledges the kind hospitality of Université de Tours during a sabbatical, funded by the CAPES/COFECUB Program (grant no.~88887.879175/2023-00), and of the Polytechnic University of Bari, funded by the INdAM--GNAMPA Visiting Program 2026 (CUP:~E53C24001950001). He also thanks Ricardo Freire for fruitful discussions regarding the topic during a visit to the CMM at the Universidad de Chile.

N.~De Nitti is a member of the Gruppo Nazionale per l'Analisi Matematica, la Probabilità e le loro Applicazioni (GNAMPA) of the Istituto Nazionale di Alta Matematica (INdAM) and has received support from the INdAM--GNAMPA Project 2026 \textit{Modelli Non-locali in Fluidodinamica, Traffico ed Elasticità} (CUP:~E53C25002010001).

\medskip
\paragraph{Use of AI tools.} OpenAI's GPT-6~Astra and Anthropic's Claude
Fable~5.1 were used to assist in the revision of the manuscript. The authors developed and verified all results and proofs and take full
responsibility for the content of this paper.

\printbibliography

\vfill 

\end{document}